\documentclass[11pt]{amsart}

\usepackage{amsthm,amsmath,natbib}
\RequirePackage[colorlinks,citecolor=blue,urlcolor=blue]{hyperref}

\usepackage{bbm}
\usepackage{booktabs}
\usepackage{geometry}                
\usepackage{graphicx}
\usepackage{amssymb}
\usepackage{setspace,multirow}
\newtheorem{theorem}{\rm T{\footnotesize{HEOREM}}}

\newtheorem{proposition}{\rm P{\footnotesize{ROPOSITION}}}
\newtheorem{condition}{Condition}
\newtheorem{corollary}{\rm C{\footnotesize{OROLLARY}}}
\newtheorem{lemma}{\rm L{\footnotesize{EMMA}}}

\DeclareMathOperator*{\argmin}{arg\,min}

\newcommand{\ba}{\mbox{\bf a}}
\newcommand{\bb}{\mbox{\bf b}}

\newcommand{\be}{\mbox{\bf e}}

\newcommand{\bh}{\mbox{\bf h}}

\newcommand{\bv}{\mbox{\bf v}}

\newcommand{\bx}{\mbox{\bf x}}

\newcommand{\bA}{\mbox{\bf A}}

\newcommand{\bE}{\mbox{\bf E}}

\newcommand{\bM}{\mbox{\bf M}}

\newcommand{\bR}{\mbox{\bf R}}

\newcommand{\bV}{\mbox{\bf V}}

\newcommand{\bX}{\mbox{\bf X}}

\newcommand{\bbeta}{\mbox{\boldmath $\beta$}}

\newcommand{\bgamma}{\mbox{\boldmath $\gamma$}}

\newcommand{\bPsi}{\mbox{\boldmath $\Psi$}}

\newcommand{\bomega}{\mbox{\boldmath $\omega$}}

\newcommand{\bDelta}{\mbox{\boldmath ${\boldsymbol\Delta}$}}

\usepackage{lipsum}
\makeatletter
\g@addto@macro{\endabstract}{\@setabstract}
\newcommand{\authorfootnotes}{\renewcommand\thefootnote{\@fnsymbol\c@footnote}}%
\makeatother

\begin{document}



\begin{center}
  \Large  {\sc
  Structured Estimation  in Nonparameteric  Cox Model \par \bigskip}

  \normalsize
  \authorfootnotes
  Jelena Bradic \textsuperscript{$\star$}
   and
  Rui Song \textsuperscript{$\dagger$} \par \bigskip

  \textsuperscript{$\star$}Department of Mathematics, University of California, San Diego \par
  \textsuperscript{$\dagger$}Department of Statistics,  North Carolina  State University \par \bigskip

\end{center}

\begin{abstract}
  
In this paper, we study   high dimensional non-parametric  Cox proportional hazard's model in a non-asymptotic setting.
We  study the finite sample oracle $l_2$  bounds for  a general  class  of group penalties that allow possible hierarchical and overlapping structures. We approximate the log partial likelihood with a quadratic functional and use truncation arguments to reduce the error.
  Unlike the existing literature, we exemplify differences between bounded and possibly unbounded non-parametric covariate effects.  In particular,   we show that  bounded  effects can lead to prediction bounds similar to the  simple linear models, whereas unbounded effects can lead to larger prediction bounds. In both situations we do not assume that  the  true parameter is necessarily sparse.   Lastly, we present new  theoretical results for hierarchical and smoothed estimation in the non-parametric Cox model, as two examples of the proposed general framework.
   \end{abstract}

\section {Introduction}  
Prediction of an instantaneous rate of
occurrence of   events
 when covariates are   high dimensional  plays a critical role in contemporary genetics studies underlying the causes of many incurable diseases. The challenge of high-dimensionality and arrival of high throughput bioinformatics data give rise to the surge of interest in the statistical literature. Most of the research in this area has concentrated on penalized methods for  Cox hazards models with parametric  \citep{T97,FL01}  or additive form  \citep{KN12,GG11}. 

Among the first theoretical work  on the Cox model with right censored data that allows $p \gg n$ is \cite{BFJ11}, where   the authors  documented good asymptotic variable selection properties  of the LASSO and SCAD penalty.   However, as with all asymptotic statistical properties, it is important to assess its relevance to the finite sample regime. 
\cite{HZ12} show non-asymptotic    oracle   estimation error   bounds for the Cox model and  with the  LASSO penalty.  Mentioned results are  derived under the assumption that the true model is exact sparse, which may be difficult to meet in practice.  One of the first work that addresses non-asymptotic oracle bounds for the nonparametric Cox model is \cite{L00} where author exemplifies importance  of   asymptotic   versus non-asymptotic theory.
\cite{KN12} analyze non-asymptotic oracle prediction bounds for the additive Cox model with fixed design, whereas \cite{GG11} analyze non-asymptotic oracle estimation bounds  for the additive hazards models.  \cite{L12} shows non-asymptotic oracle bounds for the baseline hazards function in the additive Cox model with fixed design. These work do not  assume exact sparsity and address oracle bounds for the  weighted Lasso penalty.

Censored high dimensional data are often collected from clinical studies where genomic formations are highly  complex with a large number of possible  interactions.   
Despite of its importance, the  structured sparsity  is rarely studied in the context of  censored observations. \cite {WNZZ09} give an interesting empirical study in cases of $p \leq n$.   This article focuses on censored data with   structured  (group or hierarchical)   sparsity  where covariates size is allowed to depend on dimensionality $p$ and where the sparsity is assumed to only hold approximately.    In particular,  our results easily  extend to    situations  where group LASSO,  hierarchical LASSO, group Ridge, Elastic Net  and  block $l_1/l_{\infty}$ penalty  are employed.

 To simplify the problem with essential ingredients, we assume the following nonparametric Cox model where we   let $T$  denote   the individual event time. 
  Conditional on the high dimensional covariate  $\bx$, the hazard function is modeled as 
\begin{equation}\label{eq:hazard}
\lambda(t|\bx)= \lim_{h \to 0+} \frac{1}{n} P \left( T \in[t,t+h)|T\geq t, \bx \right) =\lambda_0(t) \exp\{ g(\bx)\},
\end{equation}
for a baseline hazard function $\lambda_0(t)$ and the relative risk function $g(\bx)$.  This class of models belongs to the proportional hazards models with the relative risk function taken to be the exponential function. 
In order to estimate $g$ in high dimensional setting where the dimensionality of the covariate $\bx$ is much greater than the number of samples ($p \gg n$), it is commonly assumed that function $g(\bx)$ exhibits some form of  sparsity. In the Cox model it is typically assumed that $g(\bx)={\bbeta^*}^T\bx$ for an exactly  sparse vector $\bbeta^*$. However, this assumption  does not admit good interpretation when   $g(\bx)$  is not linear in nature. Instead, we  divide the function $g(\bx)$   into an additive   and non-additive component, where we assume that  the   additive component alone can be well approximated through sparse structured additive function.  This assumption is reasonable and more general than previous assumptions.

Exploring the structured sparsity   based on the right censored data is different from the conventional set-up hence very challenging. We approach this problem by introducing a very general class of penalty functions that allows  overlapping and group structures and  is able to control the model misspeciations mentioned above. Instead of estimating the nonparametric function $g$, we hope to find an approximate solution $\hat \bbeta$ to this problem   within a given convex subset  $\mathcal{B}$. Such a solution should minimize   the squared estimation error $ \| f_{ {\boldsymbol\beta}} - g \|^2$, where $ f_{ {\boldsymbol\beta}}$ is a linear approximation of $g$ indexed by a parameter $ { {\boldsymbol\beta}}$.

The contributions of our paper are  three folds. First, we establish two new   oracle inequalities (OI) for the high-dimensional nonparametric Cox model \eqref{eq:hazard} that explicitly bound  the squared estimation error and allow  deviations from the exact sparsity. 
Second, we show two kinds of  finite sample  sandwich bounds for the partial likelihood.  These bounds reflect  how the geometry of the partial likelihood in low dimensional space differs from that in high-dimensional space.
Third, we  show new bounds for hierarchical and smooth selection in the context of additive Cox models. In particular we discuss the complete CAP family as introduced in \cite{ZRY09} and penalties based both on sparsity and smoothness constraints, such is the elastic net penalty \citep{ZH05}, for example.

The rest of the paper is organized as follows.
In Section \ref{sec:SCGS}, we define a new class of group penalty functions and present our theoretical results in Section \ref{main-result}. 
New bounds on the distance between the least squares and the partial likelihood loss function are presented in Section \ref{sec:lan}.   Section \ref{sec:examples} is left for two examples of Hierarchical Lasso and Elastic Net penalties. 

We use the following notations.  A $pd$ dimensional vector $\bx$ is represented as  $\bx= (\mathbf x_{1}^T, \cdots, \bx_p^T)^T$ with $\bx_j =(x_{j1},\cdots, x_{jp})^T$. For a $d$ dimensional vector $\bx$, norm $\| \bx\|_{\gamma_j} = (\sum_{k=1}^d |x_{k}|^{\gamma_j})^{1/{\gamma_j}}$,  with $\gamma_j \geq 1$ and  $\gamma_j^*$ its  H\"{o}elder conjugate and such that $1/\gamma_j + 1/\gamma_j^*=1$.    Euclidian functional norm, $\| \cdot\| ^2$, is defined as
$
\| f (\bX)  \| ^2 =  \frac{1}{n} \sum_{i=1}^n  f ^2(\bX_i).
$ Throughout the paper, we denote with $\otimes$ the outer product between vectors with $\bx^{\otimes 2}$
 denoting  $\bx  \bx^T$, for any  vector $\bx$ .

\section{Convex Group Selection}\label{sec:SCGS}

Let $T$ denote the event time, $D$ denote the censoring time, and $\bX=(\bX_{1}^T,\dots,\bX_{p}^T)^T$ denote the $pd$-dimensional covariate vector where $\bX_{j} = (X_{j1},\dots,X_{jd})$. Define $Z=\min(T,D)$ and $\delta= \mathbbm{1}\{T \leq D\}$ as the observed event time and censoring indicator, respectively. We consider an i.i.d sample $\{(\bX_i, Z_i, \delta_i): i=1,\dots,n\}$ from the population $(\bX,Z, \delta)$, where $\bX_i = (\bX_{i1}^T, \dots, \bX_{ip}^T)^T$ and $\bX_{ij} = ( X_{ij,1},\dots, X_{ij,d})^T$. Let the event time, $T$, and the censoring time, $D$, be independent  conditional on the covariates. Assume that $X_i \in [a,b]$, for some constants $a$ and $b$.
Let us denote with $t_1<\cdots<t_N$   the ordered failure times and with $\mathcal{R}_q=\{ i\in\{1,\dots,n\}: Z_i \geq t_q\}$   the at risk set at each failure time.
 We define counting processes  $N_i(t)=1\{Z_i \leq t, \delta_i=1 \}$,  $ \bar{N}(t) = n^{-1}\sum_{i=1}^n N_i(t)$  and predictable processes $Y_i=1\{Z_i \geq t\} \in [0,1]$. It   holds that $dN_i(t)=dM_i(t) + d \Lambda_i(t)$ with  a martingale sequence $M_i$ and a compensator 
$$d \Lambda_i(t)=\lambda_0(t) \exp\{g(\bX_i)\}dY_i(t), $$ where $g(x)$  is the unknown function of interest. 
 Moreover, we use $\Lambda_0(\tau)=\int_0^\tau \lambda_0(t) dt$ to denote the integrated baseline function.
 
 To approximate $g(\bX)$, we define a collection of univariate   functions  $\{f_1(x),\cdots, f_p(x)\}$. We also define  a collection of dictionary functions $ \{ \Psi_1(x), \cdots, \Psi_d(x)\}$ with examples including wavelets, splines, step functions,  frames etc. 
We aim to approximate $g(\bX)$ as a linear combination of univariate functions $f_j$, each of which we  approximate  with a  linear combination of dictionary functions $\Psi_k(x)$. In more details, we approximate  $g(\bX)$ with
$$
  f_{\mathbf b}(\bX_i) = \sum_{j=1}^p f_j(X_{ij}) = \sum_{j=1}^p  \sum_{k=1}^{d} b_{jk} \Psi_{k}(X_{ij}) = \bb^T \bPsi(\bX_i) ,
  $$
  where   $ \bPsi(\bX_i) =(\bPsi(X_{i1})^T,\cdots, \bPsi(X_{ip})^T)^T$ with $ \bPsi(X_{ij})=(\Psi_1(X_{ij}), \cdots, \Psi_d(X_{ij}))^T$.
The candidate functions  are known a priori with $|\Psi_k(x) |\leq C<\infty$,  but need not be orthogonal. Note that we do not make assumptions on the number of candidate functions $d$ or $p$ and we allow both to grow and be much larger than $n$.
Let $\tau$ denote the end of the study time, we define the empirical risk function
$
{\mathcal R}_n(\bb) = - \mathcal{L}_n ( \bb, \tau)
$
and $\mathcal{L}_n (\bb,\tau)$ denotes
the  log  partial  likelihood associated to the additive component using the counting process notation:
\begin{equation}\label{eq:5}
\mathcal{L}_n (\bb,\tau)=  \frac{1}{n} \sum_{i=1}^n  \int_0^\tau   f_{\mathbf b}(\bX_i)   dN_i(t)-  \int_0^\tau  \log \mathcal{S}^{(0)}_n(\bb,t) d \bar{N}(t),
\end{equation}
with 
 $$
 \mathcal{S}^{(l)}_n(\bb,t)=\frac{1}{n}\sum_{i=1}^nY_i(t)  \bPsi^{\otimes l}(\bX_i)\exp\{f_{\mathbf b}(\bX_i) \}, \qquad l=0,1,2.
$$
We denote  population equivalents   of $ \mathcal{S}^{(l)}_n(\bb,t)$ with $ {s}^{(l)} (\bb ,t)= E_{Y,X} \mathcal{S}^{(l)}_n(\bb ,t)$. We also define 
$\mathcal{S}^{(0)}_n(g,t)=\frac{1}{n}\sum_{i=1}^nY_i(t)  \exp\{g(\bX_i) \}$ to denote censored empirical average of the unknown hazard function.   
    Later on we denote  $\mathcal{L}_n ( \bb, \tau)$ as $\mathcal{L}_n ( \bb)$ for simplicity.
    
  %

 With that in mind,    we  consider a class of estimators $\widehat \bbeta$ that solve the following penalized problem
\begin{equation}\label{eq:min1}
\widehat{\bbeta}=\argmin _{ \mathbf b \in \mathbb{R}^{pd}} \left \{ - \mathcal{L}_n (\bb) +  \lambda_n P(\bb)\right\},
\end{equation}
where we define the group penalty function (GPF), $P(\bb)$, as
\begin{equation}\label{eq:pen1}
P( \bb) = \sum_{j=1}^p \ d^{1/\gamma_j^*}\cdot \rho \left( \| \bb_j\|_{\gamma_j} \right),
\end{equation}
with a convex function $\rho$.
  The scaling, $d^{1/\gamma_j^*}$, ensures that the penalty term and the number of parameters within each group are of the same order.

We fix some vector $\bbeta^* \in \mathcal{B}$ such that $\bbeta^* = (\bbeta_1^*, \bbeta_2^*, \ldots, \bbeta_p^*)^T$,
$\bbeta^*_j  \neq   0$, for $j\in \mathcal{M}_*, ~\|\bbeta^*_j\|_{\gamma_j}  = 0 , j \in \mathcal{M}_*^c $. Set    $\mathcal{M}_*$ is any subset of $\{1,\dots, p\}$ that has the most $s$ elements, i.e., such that $ |\mathcal{M}_* |\leq s. $  Such a vector posses structured or grouped sparsity. \cite{KN12} work with all sparse vectors $\bbeta^*$ such that $\| \bb  - \bbeta^*\|_1\leq M$ for any vector $\bb$ in the parameter space and a constant $M$ independent of dimensionality $p$. Since we consider   parameter spaces that expand with the dimension, we choose vector $\bbeta^*$  among all structured-sparse vectors, such that the oracle estimator $f_{\boldsymbol\beta^*}$
  is the closest, in  the Euclidean distance, from the unknown function $g$,
   i.e. such that $ \| f_{\boldsymbol\beta^*} - g\|^2 = \min_{\mathbf b} \| f_{\mathbf b} - g\|^2$. 
    Notice that $f_{\boldsymbol\beta^*}$ is an oracle estimator  as function $g$ is unknown to us and that $\| f_{\boldsymbol\beta^*} -g \|^2=0$  if and only if  $f_{\boldsymbol\beta^*}=g$ almost surely. This does not  impose an additive structure for the true function $g$  itself.

  The following property of the introduced GPF is important in establishing finite sample bounds. We leave the  proof to the Appendix.
  
  \begin{lemma}\label{lemma:min1}
Let $\bv  =(\bv_{ 1}^T,\cdots,\bv_{ p}^T)^T \in \mathbb{R}^{pd}$, with $\bv_{n,j} \in \mathbb{R}^d$. 
Let  
$\mathcal{E}_{n,j}=\left\{\|\bv_{ j}  \|_{\gamma_j^*}\leq \lambda_n   d^{1/\gamma_j^*}\rho'(0+) \right \}$. Then,  if all the   events $\mathcal{E}_{n,j}$ hold with $j =1,\dots,p$, we have that  the  GPF family \eqref{eq:pen1}  with convex functions $\rho$  satisfies
\begin{equation}\label{eq:temp11}
{\bbeta^*}^T \bv   =\min_{ \mathbf x \in \mathbb{R}^{pd}} \left\{\lambda_n P(\bx)-(\bx-\bbeta^*)^T \bv  \right\} ~\mbox{and}
~ 
|{\bbeta^*}^T \bv  |=\min_{ \mathbf x \in \mathbb{R}^{pd}} \left\{\lambda_n P(\bx)-| (\bx-\bbeta^*)^T \bv | \right\}.
\end{equation}
\qed
\end{lemma}

The GPF
includes a wide variety of  grouping structures: $\rho$
determines how groups relate to one another, while $\{l_{\gamma_j}\}_{j=1}^p$ norms dictate the relationship of the coefficients among each group, $j$.
For $\rho=l_1$ and any $\gamma_j$, the penalty function reduces to the CAP family
 of \cite{ZRY09};
  for $\rho=l_1$ and $\gamma_j=2$, it becomes the  group Lasso penalty
 of \cite{YL06};  for $\rho=l_1$ and $\gamma_j=\infty$ it reduces to the block $l_1/l_{\infty}$ penalty
of   \cite{NW09}. The problem can be reparametrized to include a variety of  scaling factors in the penalty function.  For example,  $\rho \left( \| \bb_j^T \bR_j\|_{\gamma_j} \right)$, with proper weights $\bR_j$, or $\rho \left( \| \bR_j \bb_j\|_{\gamma_j}  + \sqrt{\bb_j^T \bM_j \bb_j} \right)$, with smoothing matrix  $\{\bM_j\}_{kl} = \int \Psi_k^{''}(x_{j}) \Psi_l^{''}(x_{j}) d x_j$ \citep{MGB09}. In section \ref{sec:examples}, we discuss these  cases in detail.

\section{Main Results}\label{main-result}

In this section, we present the main results and establish the non-asymptotic oracle inequalities  of  $\widehat \bbeta$ in terms of the  $l_2$ prediction error.  Our results differ from the previous literature in terms of the penalty function  and the measure of prediction error.
We present non-asymptotic prediction properties that allow the number of covariates  to depend and $n$ while allowing  complicated group structures in the model.  
Most of existing theoretical derivations in literature are based on   the assumption of bounded covariates, defined in Condition \ref{cond:bounded}.

\begin{condition}\label{cond:bounded}
 There exists a $M_p < \infty$ such that $\sup_{\boldsymbol b \in \mathcal{B}} \exp\{ f_{\mathbf b}(\bX_i)\} \leq M_p$.
\end{condition}
Such a  condition is often  assumed in studies where the dimension of the covariates is considered as fixed, but  should be carefully addressed in high dimensional settings where   $p \geq n$.  For example in cases where $\mathcal{B}$ is a compact $p$-dimensional ball of radius $r$ and $\bX_i$s are i.i.d. standard gaussian,  then $\log M_p =  r^2  \sqrt{\log p /n}  $  is unbounded for all $r \geq n^{1/4}$.
 Moreover, most of finite sample studies rest on a fixed design setup, a condition rarely satisfied in large genomics studies with the presence of censoring. 
  Most of this paper is dedicated to develop theory that allows deviations from such a Condition \ref{cond:bounded} in a random design setting. We present two finite sample results, where the first is rested  on Condition \ref{cond:bounded} (Theorem 1) whereas, the second  isn't  requiring such   condition (Theorem 2).   Moreover, our work does not require the quadratic structure to be present in the log partial likelihood. 
  In cases where  $p \geq n$ there is very little evidence that such a structure exists and is nontrivial.   
  For the case of a Cox model,   \cite{J89} showed that the partial likelihood is bounded away from zero if and only if the covariate vectors $\bX_i$ span the entire $p$-dimensional space, a condition that  fails if $p \geq n$.
Minimum of the log partial likelihood in the neighborhood needed to derived the theoretical properties may be too close to zero when $p$ gets large. This   indicates that prediction measure based on  Kullback-Leibler divergence for example, may falsely report small upper bounds.  Hence, we present error bounds of $l_2$ type as a usual benchmark from linear regression.

To present the results we  define 
$$
\bE_n( \bb,t)=S_{n}^{(1)}( \bb,t)/S_n^{(0)}( \bb,t),
\bV_n( \bb,t) = {S_{n}^{(2)}( \bb,t)}/{S_n^{(0)}( \bb,t)}
-\Bigl(\bE_n( \bb,t)\Bigr)^{\otimes 2},
$$
and with them the 
  gradient and the Hessian of the log partial likelihood   
$$
 \bigtriangledown \mathcal{L}_n ( \bb)   = - n^{-1}\sum_{i=1}^n \int_0 ^\tau \left( \bE_n( \bb,t) -   \bPsi (\bX_i) \right) dN_i(t), \ 
 - \bigtriangledown^2 \mathcal{L}_n ( \bb)  = n^{-1}\sum_{i=1}^n \int_0^\tau  \bV_n( \bb,t) d N_i(t).
$$
Unlike the traditional Cox models,  $ \bigtriangledown \mathcal{L}_n ( \bbeta^*)$ is not a martingale as the compensator $d \Lambda_i(t)$ does not vanish.
The following condition replaces  the classical conditions \citep{FH05} used in the asymptotic analysis of the estimation properties of the Cox model, such as those presented  in Condition 2 of \cite{BFJ11}.

\begin{condition}\label{cond:survival}
The nonparametric function of interest has bounded expectation i.e. 
  $E \exp\{g(\bX_i)\} < \infty.$
Moreover, the process $Y(t)$ is left continuous with right hand limits and such that 
  $D:=P(Y(\tau)=1 ) >0$ and 
  $\Lambda_0(\tau) < \infty$.


\end{condition}
 


Before we state the main oracle inequality  we provide concentration of measure  for the gradient of the log partial likelihood at the sparse vector $\bbeta^*$.
To that end, we need a preliminary result of Lemma \ref{lemma:temp1} providing concentration of measure for the vector $\bE_n(\bbeta^*,t)$, representing the expectation of a covariate vector $\bPsi(\bX_i)$ if an individual $i$ was selected with probability 
\[
Y_i(t) \exp\{{\bbeta^*}^T \bPsi(\bX_i)\} / \mathcal{S}_n^{(0)}(\bbeta^*,t),
\]
i.e., if an individual $i$ is selected  at time $t$  with a probability proportional to his or her intensity.

\begin{lemma}\label{lemma:temp1}
If Condition \ref{cond:survival} is satisfied, then there exists a constant $W>0$ independent of  $p,n,$ and $d$, such that  for every sequence of positive numbers $r_n$,
\begin{equation}
P \left( \sup_{0\leq t \leq \tau} \left\| \bE_n(\bbeta^*,t)-    \frac{s^{(1)} (\bbeta^*,t)}{s^{0} (\bbeta^*,t)} \right\|_\infty  \geq  c r_n + \sqrt{\frac{\log 2d }{n u^2}}  
  \right) \leq  \frac{3}{8ed}  W^2  e^{-n r_n^2D^2/u^2 e^{2m^*C}}  +  e^{-  n D^2 /2 },
  \end{equation}
  for $ \log u =  \| \bbeta^* \|_1 $ and
  $c=   1+2    \exp\{m^* C - \log D+ C \log u \}    $,
   with $m^*$ being the minimal signal strength defined as
  $m^* = \min \{ \| \bbeta^*_j\|_{\gamma_j}: j\in \mathcal{M}_* \}$.
  \qed
\end{lemma}

Careful inspection of the proof of 
Lemma \ref{lemma:temp1}, reveals that   both Condition 2(i) and Condition 2(iv) of \cite{BFJ11} now hold  with high probability, hence is of importance for the analysis of the Cox related models. 

 The next result gives tail probabilities that will be used to control the approximation error. 
 They both  depend on the GPF and require nontrivial proofs. 
Our theoretical derivations are further complicated due to the lack of martingale structure in the score vector $ \bigtriangledown \mathcal{L}_n ( \bbeta^*) $. Denote $\bh_n(\bbeta^*) =- n^{-1}\sum_{i=1}^n \int_0 ^\tau \left( \bE_n( \bb,t) -   \bPsi (\bX_i) \right) dM_i(t)$. We have the following result:

\begin{lemma}\label{lemma:eventbounds}
If Condition \ref{cond:survival} is satisfied, then for $ M = 1/(\tau \lambda_0(\tau) \Lambda_0(\tau ) C)$ and a  sequence of positive numbers $\lambda_n$ and all 
    $j =1,\cdots, p$,
\begin{eqnarray}
P \left( \lambda_0(\tau)\left|    \int_0^\tau S_n^{(0)}(g,t)  dt \right|    \|  \bPsi(X_{ij})   \|_{\gamma_j^*}\geq \lambda_n  d^{1/\gamma_j^*}\rho'(0+)  \right)\\
\nonumber
 \leq e^{ - \frac{n^2 M^2 \lambda_n^2\rho'(0+)^2}{2 \theta^2 + 2 M\sqrt{n} \lambda_n \rho'(0+)y/3}} + P(\max_{1 \leq i \leq n} \exp\{g(\bX_i)\} >y )
\end{eqnarray}
for   a 
  truncation value $y$ such that 
$
\theta^2 \geq \sum_{i=1}^n E{ \exp\{2g(\bX_i) \}1\{ \exp\{2g(\bX_i) \} \leq y\}}.
$
Moreover, there exists a constant $W>0$ independent of  $p,n,$ and $d$, such that   
for  
\begin{equation}\label{eq:C}
{C}_{\lambda_n,n,p,d}=\min \left\{ \frac{ C \lambda_n \rho'(0+)}{2 \lambda_0(\tau)},  \frac{  C \lambda_0(\tau) D^2 d^{2/\gamma_j^*} \log d}{u^4 e^{2 m^* C}}, \frac{D^2 }{2n},  \frac{ M^2 \lambda_n^2\rho'(0+)^2}{2 \theta^2 + 2 M\sqrt{n} \lambda_n \rho'(0+)y/3} \right\},
\end{equation}
we have 
\begin{eqnarray}
P \left( \|\bh_{n,{j}}(\bbeta^*)\|_{\infty} \geq  \lambda_n d^{1/\gamma_j^*} \rho'(0+)  \right)  \hskip 190pt\\ \nonumber
 \leq   2pd  \left (\max\biggl\{   \left(3+ \frac{3W^2 }{8ed}\right)e^{-n^2 {C}_{\lambda_n,n,p,d}} ,    e^ {   -n \frac{ \lambda_n^2 d^{2/\gamma_j^*} {\rho'}^2(0+)}{16 c_1^2 C^2  u^2} } \biggl\}+P(\max_{1 \leq i \leq n} \exp\{g(\bX_i)\} >y)  \right).  
\end{eqnarray}
 
\end{lemma}

 For the clarity of exposition the proof is relegated to the Appendix B.

To establish the sparse oracle inequality, we  assume a restricted eigenvalue condition similar to \cite{BRT09}. The {\bf Restricted Eigenvalue assumption, RE($\mu,s,\rho,\bgamma$)} is introduced as follows.

\begin{condition}\label{cond:re}
There exists a positive number  $\zeta=\zeta(s)>0$ such that
  \begin{eqnarray}\label{eq:RE}
\min_{\bx \in \mathbb{C}_{\mu,\rho}, \boldsymbol{\boldsymbol x} \neq 0 }  - \frac{{\mathbf x}^T \bigtriangledown^2 \mathcal{L}_n(\bbeta^* ) {\bx}}{\sum_{j\in\mathcal{M}_*} \rho(\| {\boldsymbol x}_j\|_{\gamma_j})^2}
\geq \zeta^2,
    \end{eqnarray}
    where
$  \mathbb{C}_{\mu,\rho}=\left\{ \bb \in \mathbb{R}^{pd}:  P(\bb_{\mathcal{M}_*^c})  \leq \mu  P(\bb_{\mathcal{M}_*}) \right\},
 $
  for   $\mathcal{M}_* = \{ j \in \{1,\dots, p\}: \|\bbeta_j^*\|_{\gamma_j}\neq0 \}$, 
  $|\{\mathcal{M}_*\} | \leq s$.
  \end{condition}

   The set $\mathbb{C}_{\mu,\rho}$ consists of all vectors that have support similar to the sparse vector, $\bbeta^*$. In particular, vectors with   more than $s$ non-zero elements also belong to the cone $\mathbb{C}_{\mu,\rho}$. We only require  that  their components positioned outside of 
   $\mathcal M_*$  are  smaller in size  than their  components positioned inside  $\mathcal M_*$.
  For example,  if  $\rho=l_1$ the set $\mathbb{C}_{\mu,\rho}$ is a cone formed by all vectors $\bb$ satisfying  $\| \bb_{\mathcal M _*^c}\|_1 \leq \mu \| \bb_{\mathcal{M}_*} \|_1$ as defined in \cite{BRT09}. For $\rho=l_1$ and $\gamma_j=2$, $\mathbb{C}_{\mu,\rho}$ is  the cone formed by all vectors $\bb$ satisfying  $\| \bb_{\mathcal M _*^c}\|_2 \leq \mu \| \bb_{\mathcal{M}_*} \|_2$ as defined in \cite{LPTG10}. Its geometry  changes with the penalty function, $\rho$, and the chosen $\gamma_j's$. Thus, we use the notation {\bf RE}($\mu,s, \rho,\bgamma$) to describe its dependence on the sparsity size, $s$, and the choice  $\boldsymbol\gamma=(\gamma_1,\cdots, \gamma_p)^T$, the vector  of norms used to describe the ``smoothness" of each $f_j$. 
 Further discussion of this condition is relegated to the Appendix \ref{sec:RE}.

Let us introduce two constants 
$0 \leq \upsilon_1 \leq 1$ and $0 \leq \upsilon_2 \leq 1$ satisfying
 \begin{equation}\label{eq:upsilons}
\upsilon_1\exp\{-2 C \upsilon_1\}  \leq 16 \lambda_n^2   \rho'(0+)\frac{ {\bar d} }{\zeta^2} ,  \upsilon_2  \exp\{ -2 C \upsilon_2\}  - 4\lambda_n  \frac{\bar d }{ \zeta^{2} {\rho'}^{2}(0+)} \sqrt{\upsilon_2} \leq 16 \lambda_n^2  \frac{{\bar d}^{3/2}}{   \zeta^{3} {\rho'}^{3/2}(0+) },
 \end{equation}
 where $\bar d = \sum_{j \in \mathcal{M}_*} d^{2/\gamma_j^*}$.

     \begin{lemma}\label{lem:templem}
    Let $\hat \bbeta$ be defined as in \eqref{eq:min1}  with  penalty function  GPF  defined in \eqref{eq:pen1}.   Let  Condition \ref{cond:survival} and  assumption {\bf RE}($7,s, \rho,\bgamma$) hold with $\zeta=\zeta(s)$.  Then, with probability $1-\delta$, for $\delta$ in \eqref{eq:delta} and   all $\mathbf b \in \mathbb{R}^{pd}$          $$
 2 \lambda_n  \sum_{j \in \mathcal{M}_*} d^{1/\gamma_j^*}    \rho(\| \widehat \bbeta_j - \bb_j\|_{\gamma_j})
 \leq
  64\lambda_n^2   \frac{ {\bar d} }{\zeta^2}  \exp\{ 2 C \upsilon_1\}    +32 \lambda_n^2   \frac{ {\bar d} }{\zeta^2}  \exp\{ 2 C \upsilon_2 \}   ,
$$
for $0 \leq \upsilon_1,\upsilon_2 \leq 1$ satisfying \eqref{eq:upsilons}.
\qed
     \end{lemma}

  We also define 
 sequence of weight vectors $\bomega(\bb)=(\omega_1(\bb),\cdots, \omega_n(\bb))^T$  as follows,
 \begin{equation}\label{eq:weights}
 \omega_i(\bb)=\sum_{q=1}^N \frac{\exp\{{\bb}^T \bPsi(\bX_i)\} 1\{i\in \mathcal{R}_q\}}{ \sum_{l \in \mathcal{R}_q} \exp\{{\bb}^T\bPsi(\bX_l) \}}.
 \end{equation}

 With these preparations we are ready to state the main result for the case of bounded covariate effects.

   \begin{theorem}\label{cor:localSOI}
Let $\hat \bbeta$ be defined as in \eqref{eq:min1} and    penalty function $P(\bb)$  defined in \eqref{eq:pen1}.   Let  Conditions \ref{cond:bounded} and  \ref{cond:survival} hold. Let assumption {\bf RE}($7,s, \rho,\bgamma$) hold with $\zeta=\zeta(s)$. Then,
 for any non-negative constant $A>0$ and  $\log u =\| \bbeta ^* \|_1 $, and 
 \[
 \lambda_n \geq \frac{8 A u n^{1/4} \lambda_0(\tau)}{d \rho'(0+)} \sqrt{\frac{\log pd}{n}},
 \]
\begin{equation}\label{eq:localSOIa}
\| f_{\widehat{{\boldsymbol \beta}}} - g\| ^2
 \leq 
  \min_{\mathbf b \in \mathbb{R}^{pd},  |\mathcal{M}_* | \leq s}  \left\{   (1+ \underline \omega^{-1} )   \| f_{ \mathbf b} -  g\|^2 +   64\lambda_n^2   \frac{ {\bar d} }{\zeta^2 \underline \omega }  \exp\{ 2 C \upsilon_1\}    +32 \lambda_n^2   \frac{ {\bar d} }{\zeta^2  \underline \omega }  \exp\{ 2 C \upsilon_2 \} \right \},
\end{equation}
 with  with  probability no less than $1-\delta$, $\delta>0$, where
\begin{eqnarray}\label{eq:delta}
 \delta &=& 2pd   \max\biggl\{   \left(3+ \frac{3W^2 }{8ed}\right)e^{-n^2 {C}_{\lambda_n,n,p,d}} ,    e^ {   -n \frac{ \lambda_n^2 d^{2/\gamma_j^*} {\rho'}^2(0+)}{16 c_1^2 C^2  u^2} }, e^{ -n^2 \frac{ M^2 \lambda_n^2\rho'(0+)^2}{2 \theta^2 + 2 M\sqrt{n} \lambda_n \rho'(0+)y/3} } \biggl\}
 \\ \nonumber
 &+&4pdP(\max_{1 \leq i \leq n} \exp\{g(\bX_i)\} >y),
\end{eqnarray}
 for 
 $\theta,y,M,u,m^*$ as in Lemma \ref{lemma:eventbounds}, and 
 $\bar d = \sum_{j \in \mathcal{M}_*} d^{2/\gamma_j^*}$,   $0 \leq \upsilon_1 \leq 1$ and $0 \leq \upsilon_2 \leq 1$ satisfying \eqref{eq:upsilons}
and 
 $$
\underline\omega=\min \left\{ \omega_i(\bbeta^* +c(\bb-\bbeta^*)):  \substack{\mathbf b \in  \mathbb{C}_{7,\rho}} , c\in(0,1),  i \in 1,\cdots,n \right\}$$
  for   $\bomega(\bb)$  in \eqref{eq:weights} and $\mathbb{C}_{7,\rho}$ in Condition \ref{cond:re}.

       \end{theorem}

     \begin{proof}[Proof of Theorem \ref{cor:localSOI}]

     This proof  requires careful analysis of the possible  model  misspecification and uses results of Propositions \ref{lem:approx} and \ref{prop:approx} stated in Section \ref{sec:lan}. 
      To that end,
we define an empirical functional norm $\| \cdot \| _{n, {\mathbf b^*} }$ for functions $f_{\mathbf b}:R^{pd} \to R$, $\bb \in R^{pd} $ with a fixed $c \in (0,1)$ and  $\mathbf b^* =c\bb+(1-c)\bbeta^*\in R^{pd}$,
\begin{equation}\label{eq:betanorm}
\| f _{\mathbf b}\| _{n,  \mathbf b^*}^2=
 \frac{1}{n} \sum_{i=1}^n \int_0^\tau Y_i(t) \omega_i (\mathbf b^*,t)   f_{\mathbf b}^2(\bX_i) d \bar N(t)
 -
 \Bigl[ \frac{1}{n}
  \sum_{i=1}^n \int_0^\tau Y_i(t) \omega_i (\mathbf b^*,t)  f_{\mathbf b}(\bX_i)  d \bar N(t)\Bigr]^2,
\end{equation}
for nonnegative weight process
\begin{equation}\label{eq:1}
\omega_i (\mathbf b^*,t) = \exp\{f_{ \mathbf b^*} (\bX_i)\}/ S_n^{(0)}( \mathbf b^*,t),
\end{equation}
and $\bar M(t) = n^{-1} \sum_{i=1}^n M_i(t)$.
This norm is connected to the curvature of the partial likelihood and is further discussed in Section  \ref{sec:lan}.
From the Taylor expansion and some algebra, we have that the following representation  holds for all $\bb$:
\begin{eqnarray}\label{eq:risk3} \nonumber
{\mathcal R}_n(\widehat { \bbeta} )-  {\mathcal R}_n( \bb)
=
  \frac{1}{2}   \| f_{\widehat{{\boldsymbol \beta}}} - f_{{\boldsymbol \beta}^*} \| _{n, \mathbf b_{\widehat { \boldsymbol\beta}}  }^2 -   \frac{1}{2}   \| f_{{\mathbf b}} - f_{{\boldsymbol \beta}^*} \| _{n, \mathbf b^*}^2
+
   (  \bb - \widehat{ \bbeta} )^T  \bh_n( \bbeta^*) +   {\boldsymbol\nu}_n(\hat{\boldsymbol\beta},\mathbf b,  g)  ,
 \end{eqnarray}
for $\bb_{\widehat { \boldsymbol\beta}}=c \widehat \bbeta + (1-c) \bbeta^*$ and $\mathbf b^*=\tilde c \bb + (1-\tilde c) \bbeta^*$  with a particular choice of  $c \in(0,1)$ and $\tilde c=\tilde c(\bb) \in(0,1)$, $\bh_n(\bbeta^*)$ as in Lemma \ref{lemma:eventbounds} and
   \begin{eqnarray}\label{eq:risk4}
 {\boldsymbol\nu}_n(\hat{\boldsymbol\beta},\mathbf b,  g) =-\frac{1}{n} \sum_{i=1}^n \int_0^\tau \lambda_0(t) Y_i(t) \exp\{g(\bX_i)\} \left( \log \mathcal{S}_n^{(0)}(\hat\bbeta,t) - \log \mathcal{S}_n^{(0)}(\mathbf b,t) \right) dt
 \\ \nonumber
 + \frac{1}{n} \sum_{i=1}^n  \int_0^\tau  \lambda_0(t) Y_i(t) \exp\{g(\bX_i)\} \left( \hat\bbeta^T \bPsi(\bX_i) - \mathbf b^T \bPsi(\bX_i) \right) dt.
 \end{eqnarray}
 where in the last expression we used  the  Doob Mayer decomposition  $d N_i =dM_i + d\Lambda_i$ with $d \Lambda_i = \lambda_0(t) Y_i(t) \exp\{g(\bX_i)\} dt$.

From the definition of the penalized estimator as the minimizer of penalized empirical risk in \eqref{eq:min1},  we obtain
$
 {\mathcal R}_n(\widehat{ \bbeta}) +  \lambda_n P(\widehat{ \bbeta}) \leq  {\mathcal R}_n({ \bb}) +  \lambda_n P({ \bb}),
  $
i.e. 
 \begin{eqnarray}
 \| f_{\widehat{{\boldsymbol \beta}}} - f_{{\boldsymbol \beta}^*} \| _{n, \mathbf b_{\widehat { \boldsymbol\beta}} }^2
 \leq  && \| f_{{\mathbf b}} - f_{{\boldsymbol \beta}^*} \| _{n, \mathbf b^*}^2
 + 2   (  \bb - \widehat{ \bbeta} )^T \bh_n( \bbeta^*) +
 2   {\boldsymbol\nu}_n(\hat{\boldsymbol\beta},\mathbf b,  g) \nonumber
 \\
 &+&
 2  \lambda_n (P({ \bb})-P(\widehat{ \bbeta})).\label{eq:risk4aa}
\end{eqnarray}

According to \eqref{eq:risk4} we decompose  ${\boldsymbol\nu}_n(  \mathbf b, \hat{\boldsymbol\beta}, g) $  in two parts, one that can be tied up with the estimation error and another that can be tied up with  the penalty term.
To that end, we observe that
\begin{eqnarray*}\label{eq:risk6}
 {\boldsymbol\nu}_n(\hat{\boldsymbol\beta},\mathbf b, g) \leq  \lambda_0(\tau)\left|    \int_0^\tau S_n^{(0)}(g,t)  dt \right|   \times  \hskip100pt\\  \left(  \sup_{t \in [0,\tau]} \left|  \log \mathcal{S}_n^{(0)}(\hat\bbeta,t) - \log \mathcal{S}_n^{(0)}(\mathbf b,t)   \right|
 +  ( \hat\bbeta  - \mathbf b)^T  \max_{1 \leq i \leq n }  \bPsi(\bX_i)   \right).
 \end{eqnarray*}
We denote 
 $\mathcal{S}^{(0)}_n(g,t)=\frac{1}{n}\sum_{i=1}^nY_i(t)  \exp\{g(\bX_i) \},$   equivalent  of $ \mathcal{S}^{(0)}_n(\bb,t)$ at the true, unknown function $g(\bx)$ and with $\lambda_0(\tau)$ denoting the value of the baseline hazard function at the end of the study time $\tau$.
Observe that $ \log \mathcal{S}_n^{(0)}(\mathbf b,t)$ is positively weighted log-sum-exp function for any value of $\mathbf b$, 
therefore it is
  Lipschitz continuous    (with constant 1 with respect to
the  $l_\infty$ norm),
 \[
\sup_{t \in [0,\tau]}  \left|  \log \mathcal{S}_n^{(0)}(\hat\bbeta,t) - \log \mathcal{S}_n^{(0)}(\mathbf b,t) \right |   \leq   \max_{1 \leq i \leq n } | \hat \bbeta^T \bPsi(\bX_i) - \bb^T \bPsi(\bX_i)|  .
 \]
  Furthermore, utilizing  Condition \ref{cond:survival}
we obtain
\begin{eqnarray*}\label{eq:risk7}
 {\boldsymbol\nu}_n(\hat{\boldsymbol\beta},\mathbf b,  g) \leq
     2 \lambda_0(\tau)\left|    \int_0^\tau S_n^{(0)}(g,t)  dt \right|     \max_{1 \leq i \leq n }  \left |  ( \hat\bbeta  - \mathbf b)^T  \bPsi(\bX_i)  \right|.
    \end{eqnarray*}

Combining with the previous result,  we get
 \begin{eqnarray}\label{eq:20}
 &&\| f_{\widehat{{\boldsymbol \beta}}} - f_{{\boldsymbol \beta}^*} \| _{n, \mathbf b_{\widehat { \boldsymbol\beta}} }^2
 \leq  \| f_{{\mathbf b}} - f_{{\boldsymbol \beta}^*} \| _{n, \mathbf b^*}^2
 + 2   (  \bb - \widehat{ \bbeta} )^T \bh_n( \bbeta^*)   \\
 &&+ 2  \lambda_0(\tau)\left|    \int_0^\tau S_n^{(0)}(g,t)  dt \right| \max_{1 \leq i \leq n }|( \hat\bbeta  - \mathbf b)^T    \bPsi(\bX_i)  |\nonumber
+
 2  \lambda_n (P({ \bb})-P(\widehat{ \bbeta})),\label{eq:risk4a}
\end{eqnarray}
for any $\bb$ and $\mathbf b^*, \bb_{\hat {\boldsymbol \beta}}$ fixed and defined as before. To show oracle inequality we need to tightly  control  the last three terms in the right hand side of the  previous inequality. The first of those is a martingale score vector at the parametric additive model, the second is measuring model misspecification whereas the third is quantifying the size  of the penalty function. 
Model misspecification are controlled by the penalty term. To that end we use the result of Lemma \ref{lemma:min1}.




  Utilizing Lemma \ref{lemma:min1} with  $\bDelta=\widehat \bbeta - \bb $ and $\bv_n =2 \bh_n(\bbeta^*)$, 
  the following holds from the first equality in \eqref{eq:temp11}
$$
{\bbeta^*}^T\bh_n(\bbeta^*)
\leq -(\bDelta - \bbeta^*)^T \bh_n(\bbeta^*)  +  \lambda_n  P(\bDelta),
$$
that is
\begin{eqnarray}\label{eq:inequality1}
4 \bDelta ^T \bh_n(\bbeta^*) \leq   \lambda_n  P(\bDelta), \qquad
\end{eqnarray}
on the event $\mathcal{E}_n$ defined as
\begin{eqnarray}\label{eq:set1}
\mathcal{E}_{n}=\bigcap_{j=1}^p\left\{2 \|\bh_{n,j}(\bbeta^*)   \|_{\gamma_j^*}\leq \lambda_n  d^{1/\gamma_j^*}\rho'(0+) \right\}. 
\end{eqnarray}
  
  Moreover, utilizing Lemma \ref{lemma:min1}  again, but now  with  $\bDelta=\widehat \bbeta - \bb $ and $\bv_n =  {4 \gamma_n}    \bPsi(\bX_i)  $, the following holds from  the second equality in \eqref{eq:temp11}
\begin{eqnarray*}
|{\bbeta^*}^T4 \gamma_n    \bPsi(\bX_i) |
&\leq&-|(\bDelta - \bbeta^*)^T4 \gamma_n    \bPsi(\bX_i) |  +  \lambda_n  P(\bDelta),
\end{eqnarray*}
on the event  $\mathcal{D}_{n,i}$ defined as
\begin{eqnarray}\label{eq:set2}
\mathcal{D}_{n,i}= \bigcap_{j=1}^p\left\{ 4 \lambda_0(\tau)\left|    \int_0^\tau S_n^{(0)}(g,t)  dt \right|    \|  \bPsi(X_{ij})   \|_{\gamma_j^*}\leq \lambda_n  d^{1/\gamma_j^*}\rho'(0+)  \right\}.
\end{eqnarray}

After rearranging the terms  and noticing that $|\bDelta^T    \bPsi(\bX_i)  | \leq  |{\bbeta^*}^T    \bPsi(\bX_i) |   +  |(\bDelta - \bbeta^*)^T    \bPsi(\bX_i) |,  $ we get
\begin{eqnarray}\label{eq:inequality2}
4 \gamma_n   |\bDelta^T  \bPsi(\bX_i)  |  \leq \lambda_n P(\bDelta),
\end{eqnarray}
  and with it that $4 \gamma_n \max_{1 \leq i \leq n }  |\bDelta^T  \bPsi(\bX_i)  |  \leq \lambda_n P(\bDelta)$, on the event $\mathcal{D}_{n,i}$.

  Therefore, combining \eqref{eq:20} with  \eqref{eq:inequality1} and  \eqref{eq:inequality2} ,  we conclude that  for all $\bb$,  conditionally on  the event
$$ \mathcal{E}_{n}  \cap \bigcap_{i=1}^n \mathcal{D}_{n,i},$$
the following inequality holds
 \begin{eqnarray*}\label{eq:risk5}
   \| f_{\widehat{{\boldsymbol \beta}}} - f_{{\boldsymbol \beta}^*} \| _{n, \mathbf b_{\widehat { \boldsymbol\beta}} }^2  \leq \| f_{ \mathbf b} - f_{{\boldsymbol \beta}^*} \|_{n,\mathbf b^*}^2 + 2 \lambda_n  \sum_{j \in \mathcal{M}_*} d^{1/\gamma_j^*} \left(   \rho(\| \widehat \bbeta_j - \bb_j\|_{\gamma_j}) + \rho(\|\bb\|_{\gamma_j})  - \rho(\| \widehat \bbeta_j \|_{\gamma_j}) \right),
  \end{eqnarray*}
  for all $\bb_{\widehat{\boldsymbol\beta}}=c \widehat \bbeta+ (1-c) \bbeta^* $ and $\mathbf b^*= \tilde c \bb + (1-\tilde c) \bbeta^*$. Let us fix $\bb_{\widehat{\boldsymbol\beta}}$ and  $\mathbf b^*$ from hereon.
  From the triangular inequality for the GPF,  we have $\rho(\| \bb_j \|_{\gamma_j}) \leq \rho(\|  \widehat \bbeta_j - \bb_j \|_{\gamma_j}) + \rho(\|  \widehat \bbeta_j  \|_{\gamma_j}) $ leading to
  \begin{equation}\label{eq:temp19}
 \| f_{\widehat{{\boldsymbol \beta}}} - f_{{\boldsymbol \beta}^*} \| _{n, \mathbf b_{\widehat { \boldsymbol\beta}} }^2
   \leq \| f_{ \mathbf b} - f_{\boldsymbol \beta^*} \|_{n,\mathbf b^*}^2 + 2\lambda_n  \sum_{j \in \mathcal{M}_*}   d^{1/\gamma_j^*}\rho(\| {\boldsymbol\Delta}_j\|_{\gamma_j})  .
     \end{equation}

 Secondly, we control the penalty term in \eqref{eq:temp19} in the Lemma \ref{lem:templem} whose   proof is presented in the Appendix E.

Utilizing further  the bound between the norms $ \| f_{ \mathbf b} - f_{\boldsymbol \beta^*} \|_{n,\mathbf b^*}^2 \leq  \| f_{ \mathbf b} - f_{\boldsymbol \beta^*} \|_{n}^2$  (proved in Proposition 4 in Section \ref{sec:lan}) in combination to \eqref{eq:risk5} , we obtain
   \begin{eqnarray*}\label{eq:temp21}
  \underline \omega  {\| f_{\hat {\boldsymbol\beta} } - f_{ {\boldsymbol\beta}^*}\|_{  }^2} \leq      \| f_{\boldsymbol\beta^*}- f_{\mathbf b}\| ^2     +  64\lambda_n^2   \frac{ {\bar d} }{\zeta^2}  \exp\{ 2 C \upsilon_1\}    +32 \lambda_n^2   \frac{ {\bar d} }{\zeta^2}  \exp\{ 2 C \upsilon_2 \},
  \end{eqnarray*}
  with $\upsilon_1, \upsilon_2$ defined above in \eqref{eq:upsilons}.
 Moreover, from the definition of the vector $\bbeta^*$ and the  triangular inequality  we have 
   \[
  \|f_{\widehat{{\boldsymbol \beta}}} - g \|^2 \leq \| f_{\widehat{{\boldsymbol \beta}}} - f_{{\boldsymbol \beta}^*} \|^2 +  \min _{\mathbf b} \| f_{{{\mathbf b}}} - g \|^2,
 \]
 which in combination to the previous inequality provides
  \[
 \|f_{\widehat{{\boldsymbol \beta}}}  - g \|^2  \leq    (1+\frac{1}{ {\underline\omega}} )  \min_{\mathbf b \in \mathcal{B}} \|   f_{\mathbf b} -g\| ^2   +  64\lambda_n^2   \frac{ {\bar d} }{\zeta^2 {\underline\omega}}  \exp\{ 2 C \upsilon_1\}    +32 \lambda_n^2   \frac{ {\bar d} }{\zeta^2 {\underline\omega}}  \exp\{ 2 C \upsilon_2 \}  .
 \]

%
%
The theorem follows easily if we bound the probability of the event $  \mathcal{E}_{n}  \cap \bigcap_{i=1}^n \mathcal{D}_{n,i}$, which   is given in Lemma \ref{lemma:eventbounds}. Hence, the proof is completed. 

 \end{proof}

 Theorem 1 establishes a new finite sample oracle inequality with possible deviations of exact sparsity.  
The first term on the right hand side of \eqref{eq:localSOIa} measures how far is the true function of interest $g(x)$ from the sparse additive approximation $f_{\boldsymbol\beta^*}$ and is only equal to zero if $g= f_{\boldsymbol\beta^*}$ almost surely. Typically, similar results appeared  in problems with fixed design \citep{KN12} or  if  one considers estimation errors related to the Kullback-Leibler divergence that are quadratic in nature   \citep{GG11,L12}. In contrast, our results hold for log partial likelihood of non-quadratic type and a class of general  random designs and general group penalty.
  The last two quantities of the RHS of \eqref{eq:localSOIa} represent the convergence rate  for the appropriate choices of $\lambda_n$. 

 Let us comment on the size of the constant $\underline{\omega}^{-1}$ appearing in the bound \eqref{eq:localSOIa}.   Each  weight, $\omega_i(\bb)$, is a sum of the conditional probabilities that observation $i$ had an event at time $t_q$, given that at least one event occurred at time $t_q$.

\begin{proposition}\label{prop:w}
Let $\eta >0$ and $c_2\in \mathbb{R}$ be constants  such that  for all $q=1,\cdots, N,$
\begin{equation}\label{eq:optcond}
\lambda_{\min}
\Biggl(
\begin{array}{cc}
\sum_{l\in \mathcal{R}_q} \bPsi^{\otimes 2}(\bX_l)   + \eta \mathbb{I}_{pd} &   \sum_{l\in \mathcal{R}_q} \bPsi(\bX_l)  \\
\sum_{l\in \mathcal{R}_q}  \bPsi^T(\bX_l)  &  c_2 - \eta b_n
\end{array}
\Biggl) = \delta^\star  ,
\end{equation}
where $\mathbb{I}_{pd}$ is a unit matrix.
Then, for all $i \in \{1,\cdots,n\}$ and $b_n >0$, the solution to the  optimization problem
\begin{equation}\label{eq:opt}
\begin{array}{cc}
\min_{ \mathbf b \in \mathbb{R}^{pd}}  \left\{ \omega_i(\bb):  \| \bb \|_2^2 \leq b_n  \right \}
\end{array}
\end{equation}
is attained and the minimum  $\omega_{\min}$ satisfies
$ \omega_{\min} \delta^\star=\sum_{q=1}^N  \min \left\{  0, \lambda_{\min} \bigl(\bPsi(\bX_i)^{\otimes 2}  \bigl) 1(i \in \mathcal{R}_q) \right\}.
$
\end{proposition}
The conditions of Proposition \ref{prop:w} are not restrictive and are easily verifiable for well posed problems.    
For $\kappa_i = \min \{  \bv \bPsi^{\otimes}(\bX_i)  \bv^T , \| \bv\|_2 \leq 1, \bv \in \mathbb{C}_{7,\rho}\}$
and by Cauchy's interlacing theorem of Hermitian matrices,    for Propositon \ref{prop:w} to hold it suffices that the random covariates $\bX_i$ satisfy $ \min _{i\in \mathbb{R}_q}  \kappa_i >0 $ .  In that case, we conclude that $\underline\omega$ satisfies
\[
\delta^\star\underline \omega \geq  \sum_{q=1}^N \min \bigl\{  0, \min _{i\in \mathbb{R}_q} \kappa_i  \bigl\}.
\]
With that, 
 for  $\kappa^q  = \min \{ \sum_{i \in \mathbb{R}_q} \bv \bPsi^{\otimes}(\bX_i)  \bv^T , \| \bv\|_2 \leq 1, \bv \in \mathbb{C}_{7,\rho}\}$
 we obtain an upper bound on the leading constant of the Theorem 1,
\begin{equation} \label{eq:varepsilon}
1+\varepsilon \leq \biggl(\sum_{q=1}^N \frac{\kappa_i }{ \kappa^q}1(i \in \mathcal{R}_q)\biggl)^{-1},
\end{equation}
with the right-hand side bounded  away from infinity almost surely. Therefore,  under Condition \ref{cond:bounded}, the proposed estimator achieves Gaussian-like oracle rates similar to those of penalized least squares   (see further discussion in Section \ref{sec:examples}).

 In the next result, we take a novel approach  and show  informative high-dimensional sparse oracle results  for possible unbounded covariate effect. First, we localize our  penalized estimator to a small elliptical neighborhood around  $\bbeta^*$.
With an appropriate choice for the tuning parameter, $\lambda_n$, the radius of the neighborhood becomes independent of the dimensionality. This is presented in the next lemma.

\begin{lemma}\label{lem:globalSOIb}
   For  $\log p \leq n$  and $s \leq \log n$, let $\hat \bbeta$ be defined as in \eqref{eq:min1} with   penalty function $P(\bb)$  defined in \eqref{eq:pen1} and $\bbeta^*$ the true sparse parameter.   Let  Condition \ref{cond:survival} and  assumption {\bf RE}($7,s, \rho,\bgamma$) hold with $\zeta=\zeta(s)$.
  If $\lambda_n$ satisfies
  \[
      \lambda_n \sum_{j\in \mathcal{M}_*} d^{1+2/\gamma_j^*}  \geq  c \zeta^4,
  \] then
 with probability $1-\delta$, for $\delta$ in \eqref{eq:delta},
\begin{eqnarray}\label{eq:cor2SOI2}
      \sum_{j=1}^p d^{1/\gamma_j^*} \| \widehat {\boldsymbol \beta}_j - \boldsymbol{\beta}^*_j \|_{\gamma_j}
      \leq 16 \sqrt{2} C e^{C+\upsilon_1} r_n ,
      \end{eqnarray}
      for $r_n= \frac{ \lambda_n }{  \zeta^2} \sum_{j\in \mathcal{M}_*} d^{2/\gamma_j^*}$ and $0 \leq \upsilon_1 \leq 1$ satisfying \eqref{eq:upsilons}.\qed
\end{lemma}

Second, using a local neighborhood structure,  we sandwich the risk with  lower and upper bounds of  quadratic form. To that end, we have the following result.

\begin{theorem}\label{cor:SOI}
For  $ \log(p d) \leq n$ let $\hat \bbeta$ be defined as in \eqref{eq:min1} and    penalty function $P(\bb)$  defined in \eqref{eq:pen1}.   Let  Condition \ref{cond:survival} and  assumption {\bf RE}($7,s, \rho,\bgamma$) hold with $\zeta=\zeta(s)$. Then,
 for non-negative constant $A>0$ and $u$ defined in Theorem 1,
\begin{equation}\label{eq:LAMBDA}
 \lambda_n \geq \frac{8 A u n^{1/4}}{d \rho'(0+)} \sqrt{\frac{\log pd}{n}} \qquad \mbox{and} \qquad \lambda_n \sum_{j\in \mathcal{M}_*} d^{1+2/\gamma_j^*}  \geq  c \zeta^4,
\end{equation}
 with probability no less than $1-\delta$, $\delta>0$ and satisfying \eqref{eq:delta},
    there exists $\varepsilon >1$ such that,
\begin{equation}\label{eq:localSOI}
\| f_{\widehat{{\boldsymbol \beta}}} - g\| _{ }^2
 \leq
 \min_{\mathbf b \in \mathbb{R}^{pd},  |\mathcal{M}_* | \leq s}  \left\{  (1+\varepsilon)   \| f_{ \mathbf b} -  g\|^2_{ }   + 64\lambda_n^2 \varepsilon  \frac{ {\bar d} }{\zeta^2}  \exp\{ 2 C \upsilon_1\}    +32 \lambda_n^2   \varepsilon\frac{ {\bar d} }{\zeta^2}  \exp\{ 2 C \upsilon_2 \} \right \},
\end{equation}
 with  $\bar d = \sum_{j \in \mathcal{M}_*} d^{2/\gamma_j^*}$, $0 \leq \upsilon_1 \leq 1$ and $0 \leq \upsilon_2 \leq 1$ satisfying \eqref{eq:upsilons}  
and $\underline{ \underline \omega } = \min_i \omega_i(\bbeta^*)$ with
$$\varepsilon =   \underline{\underline \omega}^{-1}\exp\left\{C e^C  26\frac{\lambda_n }{   \zeta^2}\sum_{j \in \mathcal{M}_*} d^{2/\gamma_j^*}\right\}.$$
 
\end{theorem}

\begin{proof}

 We proceed by first restricting the parameter space to an elliptical neighborhood that is not expanding with dimensionality $p$. Then, we
 apply  Lemma \ref{lem:templem}   (stated in the Proof of Theorem 1) and Proposition 4 (stated and proved in Section \ref{sec:lan}) to finalize the proof.

From Proposition \ref{lem:approx}, we have that
 \begin{eqnarray}\label{eqn:1a}
  \| f_{\widehat{{\boldsymbol \beta}}} - f_{{\boldsymbol \beta^*}}\|_{n,\widehat{\boldsymbol \beta }^*}^2  \geq e^{-2 a_{\widehat {\boldsymbol \beta} - \boldsymbol{\beta}^*}} \| f_{\widehat{{\boldsymbol \beta}}} - f_{{\boldsymbol \beta^*}}\|_{n,\boldsymbol\beta^*}^2 ,
 \end{eqnarray}
 with $\widehat{\boldsymbol \beta }^* = c \widehat{\boldsymbol \beta }  + (1-c)  {\boldsymbol \beta }^*$, some $c \in (0,1)$ and
$
  a_{\widehat {\boldsymbol \beta} - \boldsymbol{\beta}^*}  =  2 \max_{1 \leq i \leq n} |(\widehat {\boldsymbol \beta} - \boldsymbol{\beta}^*)^T \bPsi(\bX_i) | .
  $
The exponential term in the previous equation needs to be tightly controlled for $p \gg n$.
The proof of the theorem is then finalized by finding nontrivial bounds for the empirical norms, $\|\cdot\|_{n,\cdot}$ as defined in \eqref{eq:betanorm}, while allowing $p \gg n$. Let $p \geq n$ and $\log p \leq n$.  We establish that by bounding the appropriate norm of the error vector $\widehat{\bbeta}-\bbeta^*$ and obtaining the bound, which is log linear in dimensionality $p$. The result is summarized in the  Lemma \ref{lem:globalSOIb}, whose proof is provided in the Appendix E.

Consequently  
 we have  \begin{eqnarray*}
  a_{\widehat {\boldsymbol \beta} - \boldsymbol{\beta}^*} \leq  2  \sum_{j =1}^p  \| \widehat {\boldsymbol \beta}_j - \boldsymbol{\beta}^*_j \|_{\gamma_j}  \max_{1 \leq i \leq n}  \left( \sum_{k=1}^d (\Psi_k(X_{ij}))^{\gamma_j^*}\right)^{1/\gamma_j^*}
 \leq 32 \sqrt{2} C e^C \frac{ \lambda_n }{  \zeta^2} \sum_{j\in \mathcal{M}_*} d^{2/\gamma_j^*}.
   \end{eqnarray*}
 Remember that
   $C \geq \max_{k,i,j } | \Psi_k(X_{ij})|$. Hence,  we  have successfully localized the error  vector $\widehat{\bbeta}-\bbeta^*$ in a sparse neighborhood whose radius is  not increasing with the dimensionality $p$.

Utilizing further Lemma \ref{lem:templem} and Proposition 4 with equation \eqref{eq:temp19}, we obtain
   \begin{eqnarray*}
  \underline {\underline \omega}  {\| f_{\hat {\boldsymbol\beta} } - f_{ {\boldsymbol\beta}^*}\|_{n }^2} \leq      e^{ 32 \sqrt{2} C e^C r_n} \left\{  \min_{\mathbf b \in \mathcal{B}}  \| g- f_{\mathbf b}\| ^2     +  64\lambda_n^2   \frac{ {\bar d} }{\zeta^2}  \exp\{ 2 C \upsilon_1\}    +32 \lambda_n^2   \frac{ {\bar d} }{\zeta^2}  \exp\{ 2 C \upsilon_2 \} \right\},
  \end{eqnarray*}
  with $\upsilon_1, \upsilon_2$ defined above in Lemma \ref{lem:templem}. The proof is finalized by simple triangle inequality.
\end{proof}

Let us comment on the size of $\varepsilon$ appearing on the RHS in Theorem 2.
From the results of   Proposition \ref{prop:w},  we can easily conclude that  for $s \leq \log n$, $\varepsilon \geq 0$ and
\[
\varepsilon\leq  \exp \bigl\{   32 \sqrt{2} C e^C r_n   \bigl\},
\]
for  $r_n \to 0$ such that $r_n\zeta^2 = { \lambda_n }  \sum_{j\in \mathcal{M}_*} d^{2/\gamma_j^*}$, and the constant $C$ defined as the upper bound on the dictionary functions $\Psi$.
The difference in the rates of convergence between Theorems \ref{cor:localSOI} and \ref{cor:SOI}  reflects the dimensionality of the problem. In comparison to  Theorem \ref{cor:localSOI}, Theorem \ref{cor:SOI} differs in the presence of the exponential term of the order of  $e^{r_n}$. This additional term is coming from the complex likelihood structure  for possibly unbounded covariate effects for which we show  that   local Lipchitz constant is proportional to   $e^{r_n}$. 


Combining \eqref{eq:LAMBDA} with \eqref{eq:localSOI} we see that the oracle inequality of Theorem \ref{cor:SOI}  rely on the choice of the number of basis functions, $d$, and it  requires $   d^{-1}{\log(pd)} <  \sqrt{n} e^{-\| \boldsymbol\beta^*\|_1} /s$.   The more events we observe, the more basis functions we can choose. In the classical case of $d \sim n^{-1/2}$, the previous constraint becomes $ \log(p) \leq \frac{1}{2}\log n + \frac{n}{s }e^{-\| \boldsymbol\beta^*\|_1}$, which implies that dense problems with $p \gg n$ and $s \geq \sqrt{n}$ cannot be efficiently retrieved.

We also note that previous results do not require exact sparsity to hold, that is, they do not assume ${\boldsymbol\beta}^*$ is the true underlying parameter.  Interestingly, if we   make such an assumption, the result of Lemma \ref{lem:globalSOIb}  show that there is no effective difference in estimation between bounded and possibly unbounded covariate effects. Namely,  special case of  the result of Lemma \ref{lem:globalSOIb}  matches the result of \cite{HZ12} for lasso penalty  (consider  $\rho=l_1,\gamma_j=1$).

In summary,  the results of Theorems   \ref{cor:localSOI} and \ref{cor:SOI}  are quite general. They  cover a wide range of penalty functions with a choice of $\gamma_j$'s and are applicable to Lasso, group Lasso, group ridge, CAP penalty, elastic net and many more. Two specific  examples will be discussed in Section \ref{sec:examples}.

%
%
\section{ Sandwich Bounds for the Log Partial Likelihood}\label{sec:lan}


This section gathers some results that were crucial in obtaining Theorems \ref{cor:localSOI} and \ref{cor:SOI} . The novel  ideas of the major results are to quantify the distance between the log partial likelihood ${\mathcal R}_n( \bb)$ and the approximate quadratic expansion of the log partial likelihood. 

Without loss of generality,   ${\mathcal R}_n( \bb)$  can be written as
$
{\mathcal R}_n( \bb) =  - \mathcal{L}_n ( \bb)  + \mathcal{L}_n ( \bbeta^*)  - \mathcal{L}_n ( \bbeta^*).
$
By Taylor expansion around $ \bbeta^*$, we have that there exists a $c
\in (0,1)$ and $\mathbf b^*=c\bb+(1-c)\bbeta^*$ such that
\begin{equation}\label{eq:risk1} \nonumber
{\mathcal R}_n( \bb) =  -\left(  \bb-  \bbeta^* \right)^T \{\bigtriangledown \mathcal{L}_n ( \bbeta^*) \}  - \frac{1}{2} \left(  \bb -  \bbeta^* \right)^T \{\bigtriangledown^2 \mathcal{L}_n ( \mathbf b^*) \}  \left(  \bb -  \bbeta^* \right)
- \mathcal{L}_n ( \bbeta^*).
\end{equation}

Together with the previous Taylor expansion, the empirical risk function can be decomposed as follows.
For every $\bb$, there exists a  $c\in(0,1)$ and $\mathbf b^*=c\bb+(1-c)\bbeta^*$ such that ${\mathcal R}_n( \bb)$  admits the following quadratic representation:
\[
 {\mathcal R}_n( \bb) = -  \left(  \bb -  \bbeta^* \right)^T \{\bigtriangledown \mathcal{L}_n ( \bbeta^*) \}
 +
   \frac{1}{2}   \| f_{{\mathbf b}} - f_{{\boldsymbol \beta}^*} \| _{n, \mathbf b^*  }^2
   -
   \mathcal{L}_n ( \bbeta^*).
   \]

 Because no two counting processes, $N_i(t)$ and $N_j(t)$, jump at the same time,  the following holds:
\begin{equation} \label{eq:def1}
\| f _{\mathbf b}\| _{n,  \mathbf b^*}^2 =  \frac{1}{n}\sum_{i=1}^n \int_0^\tau Y_i(t)  \omega_i (\mathbf b^*,t)   (f_{\mathbf b}(\bX_i)-\bar f_{\mathbf b}^{*}(t) )^2    d \bar N(t),
\end{equation}
where
$
\bar f_{\mathbf b}^{*}(t) = \frac{1}{n}\sum_{i=1}^n Y_i(t)\omega_i (\mathbf b^*,t)  f_{\mathbf b}(X_i)
$
 can be understood as a process of empirical weighted averages of $ f_{\mathbf b}$.
If Condition \ref{cond:survival} is satisfied, then, there exists a $c \in (0,1)$ such that the introduced empirical norm is a proper norm. To be specific,  the norm is nonnegative, $\| f _{\mathbf b}\| _{n,  \mathbf b^*}^2=0$  for every such $\mathbf b^*$ if and only if $\bb=0$. In addition, the norm satisfies the triangular inequality that
$
  \| f_{\mathbf b_1} - f_{{\mathbf b}_2} \| _{n, {\mathbf b^*} }  \leq  \| f_{\mathbf b_{1}}  \| _{n, {\mathbf b^*} } + \| f_{{\mathbf b}_2}  \| _{n, \mathbf b^* }
$ for every $\mathbf b_{1}, \bb_2$ and  fixed $\mathbf b^*$.
%
%



As squared Eucledian norm $\|\cdot\|$   represents natural benchmark, we seek to understand the lower and upper bounds of the $\|\cdot\|_{n,\cdot}$ norm using the $l_2$ empirical norm   $\|\cdot\|$ in the next result. 

\begin{proposition}\label{lem:approx}
Let $\underline{ \underline \omega }$ be defined as in Theorem \ref{cor:SOI}.
For any vector $\bv$ define $$a_{\mathbf v}=\max_{1\leq i,q \leq n} | \bv^T [\Psi(\bX_i) - \Psi(\bX_q)] |.$$ Then, the following sandwich bound holds almost surely for every vector $\bb$ and  corresponding vector $\mathbf b^*=c\bb +(1-c)\boldsymbol \beta^*$,
 \begin{equation}\label{eq:sand1}
\underline{ \underline \omega }e^{ - 2 a_{\mathbf b- \boldsymbol\beta^*} } \ \| f_{\mathbf b} - f_{\boldsymbol\beta^*}\| ^2 \leq \| f_{\mathbf b} - f_{\boldsymbol\beta^*}\|_{n, \mathbf b^*}^2  \leq e^ { 2 a_{\mathbf b- \boldsymbol\beta^*} } \ \| f_{\mathbf b} - f_{\boldsymbol\beta^*}\| ^2  ,
\end{equation}
uniformly for every $c \in (0,1)$.
\end{proposition}

A similar result appeared independently in the recent work of \cite{HZ12} (see Lemma 4.3).
Such a result shares similarities to the self-concordant arguments of \cite{B10}, but the last arguments  do not cover  cases of $p \gg n$.


\begin{proposition}\label{prop:approx}
Let $N$ represent the number of distinct events. Then, uniformly  for every $\bb \in [-b_n, b_n]$, with $b_n >0$ satisfying the condition of Proposition \ref{prop:w},  and $\mathbf b^*=c \bb +(1-c)\bbeta^*$, with $c\in(0,1)$, the following holds almost surely:
\begin{equation} \label{eq:approx}
  n^{-1} \sum_{q=1}^N \frac{\min \left\{  0, \min _{i\in \mathbb{R}_q}  \lambda_{\min } \bigl(\bPsi(\bX_i) \bPsi^T(\bX_i)\bigl) \right\}}{\lambda_{\min}
 \bigl(\sum_{l \in \mathbb{R}_q}\bPsi(\bX_l) \bPsi^T(\bX_l)\bigl)}
\| f_{\mathbf b} \| ^2  \leq \| f _{\mathbf b}\| _{n, \mathbf b^*}^2 \leq   \| f_{\mathbf b} \| ^2    ,
\end{equation}


Moreover, if $b_n$ is bounded and $\min_{1\leq i \leq n} \lambda_{\min}
 \bigl( \bPsi(\bX_i) \bPsi^T(\bX_i)\bigl) >0$, then the left-hand bound in \eqref{eq:approx} is strictly positive almost surely.
\end{proposition}

%


Propositions \ref{prop:w}-\ref{prop:approx} are critical in establishing the main result in terms of non-trivial lower and upper bounds. We utilize Propositions  \ref{prop:w} \& \ref{prop:approx}  for low-  and  \ref{lem:approx}   for high-dimensional problems, respectively.
With the help of all four results, we are able to obtain the main results in Section \ref{main-result}.



          \section{Examples} \label{sec:examples}

In this section, we show two examples of GPFs \eqref{eq:pen1} (that allows hierarchical structures within and among groups)  and show their theoretical properties in the Cox model setup. To the best of our knowledge, similar results do not exist in the current literature. For simplicity in the presentation, the result of this section   focus on the exact sparse models with $\bbeta^*$ representing the unknown true parameter.

          \subsection{Hierarchical Selection and CAP}

Our results apply to a general class of additive  models, where the groups in the additive {\color{red} Cox} model may share some but not necessarily all features across groups. For example, the effect of one gene can be shared by many different pathways, and thus studying hierarchical gene selection is  of significant  importance. For some genes, it has already been ascertained that their over or under expression contributes to the survival rates of cancer patients.   Based on the prior information,  each $f_j$ can be approximated by  $\bb_{\Gamma_j}^T \bPsi_{\Gamma_j}$, where $\Gamma_j$ is a set of covariates that belongs to group $j$. The regularized estimator, $\widehat{\bbeta}$, is then defined as the minimizer of
         $$
          - \frac{1}{n} \sum_{i=1}^n \int_0^\tau \bigl\{\sum_{j=1}^p \bb_{\Gamma_j}^T \bPsi_{\Gamma_j}+ \log  \bigl( \frac{1}{n}\sum_{i=1}^nY_i(t)  \exp\{\sum_{j=1}^p \bb_{\Gamma_j}^T \bPsi_{\Gamma_j} \} \bigl)\bigr\} d  N_i(t)
          + \sum_{j=1}^p \lambda_{n,j} |\Gamma_j|^{1/\gamma_j^*} \|\bb_{\Gamma_j}\|_{\gamma_j},
         $$
where $|\Gamma_j|$ denotes for the cardinality of that set.          Note that this penalty includes the classical group Lasso penalty, where one would select all $\gamma_j=2$.

\begin{corollary}
Let conditions of Theorem 2 be satisfied.
Then, for some constant $A >4$ and the choice of the tuning parameters
$$ \sqrt{d} \lambda_{n,j}  \geq A    \min \left \{ \zeta^2,\sqrt{\frac{ \log(pd)}{n} |\Gamma_j|^{-2/\gamma_j^*}} \right\}, 
$$
 with the probability of at least of $1- 6\{pd\}^{1-A}$,
  \begin{equation}
\| f_{\widehat{{\boldsymbol \beta}}} - f_{\boldsymbol\beta^*}\| _{ }^2
 \leq
{\underline {\underline \omega}}^{-1} e^{   C ( r_n e^{C}  + 2\upsilon_1  )} 26 \frac{\sum_{j \in \mathcal{M}_*} \lambda_{n,j}^2 |\Gamma_j|^{2/\gamma_j^*}  }{\zeta^{2} },
\end{equation}
 for
  $r_n=   \zeta^{-2}\sum_{j \in \mathcal{M}_*} \lambda_{n,j} |\Gamma_j|^{2/\gamma_j^*} $ , and $0 \leq \upsilon_1 \leq 1$  satisfying
 \begin{equation}\label{eq:upsilons1}
\upsilon_1 e^{-2 C \upsilon_1}  \leq 16   \rho'(0+)\frac{\sum_{j \in \mathcal{M}_*} \lambda_{n,j}^2 |\Gamma_j|^{2/\gamma_j^*}  }{\zeta^2} .
 \end{equation}

\end{corollary}
The proof of this result is omitted because it is a simple modification of the results presented in the paper with $\lambda_n$ being adaptive to each group $\Gamma_j$.  To the best of our knowledge, the oracle inequality of Corollary 1 is one of first that discusses high-dimensional finite-sample properties of the whole CAP family proposed in the seminal work of \cite{ZRY09}.
In particular, the block $l_1/l_\infty$ penalty introduced in \citep{NW09}  is a member of the CAP family. \cite {NW09} present $l_\infty$
 bounds on the estimation error of block $l_1/l_\infty$ penalty in the linear models. Previous Corollary 1 provides its finite sample $l_2$ error bounds for the sparse additive Cox model with possibly overlapping groups. In more details, 
we obtain with high probability
\[
\| f_{\widehat{{\boldsymbol \beta}}_{{ l_1/l_\infty}}} - f_{\boldsymbol\beta^*}\| ^2
 \leq 26
 e^{   C ( r_n e^{C}  + 2\upsilon_1  )}  \frac{\sum_{j \in \mathcal{M}_*} \lambda_{n,j}^2 |\Gamma_j|^{2 }  }{\zeta^{2} },
\]
for  the block $l_1/l_\infty$ penalty, $0 \leq \upsilon_1 \leq 1$  satisfying
$
\upsilon_1 e^{-2 C \upsilon_1}  \leq 16   \rho'(0+) {\sum_{j \in \mathcal{M}_*} \lambda_{n,j}^2 |\Gamma_j|^{2 }  }/{\zeta^2} ,
$
and $r_n =  \zeta^{-2}\sum_{j \in \mathcal{M}_*} \lambda_{n,j} |\Gamma_j|^{2}$.   
 
Moreover, non-overlapping  groups gained significant attention with importance of   multi-task learning  \citep{LPTG10}.  Similar setup has not been investigated in models related to \eqref{eq:hazard}.   Corollary 1, provides a finite sample bound for 
  the  multi-task learning i.e.  $l_1/l_2$ penalty as follows
\[
\| f_{\widehat{{\boldsymbol \beta}}_{l_1/l_2}} - f_{\boldsymbol\beta^*}\| _{ }^2 \leq 26
  e^{   C ( r_n e^{C}  + 2\upsilon_1  )}   \frac{d \sum_{j \in \mathcal{M}_*} \lambda_{n,j}^2   }{\zeta^{2} }
\]
  with $0 \leq \upsilon_1 \leq 1$  satisfying
$
\upsilon_1 e^{-2 C \upsilon_1}  \leq 16  d \rho'(0+){\sum_{j \in \mathcal{M}_*} \lambda_{n,j}^2    }/{\zeta^2}$  and $r_n =   d  \zeta^{-2}\sum_{j \in \mathcal{M}_*} \lambda_{n,j}     .
$
If in addition Condition 1 holds, previous result  reaches the information bound of high dimensional linear models, as it matches the upper and lower bound of multi task learning as presented in \cite{LPTG10}. 

\subsection{Smooth Selection}

Throughout the previous sections, we simplified the technical details and left out the smoothing component of the penalty.  Although selection of groups of features is important,  smoothing splines  become of interest when considering non-parametric estimation. Because of  knot selections, there are potential questions of  stability of estimation.  Adding pre-described  smoothing requirements for the choice of $\Psi$  has become a standard technique for avoiding instability.   We will show that the work  of the previous sections extends to this situation with only a few adaptations. Let us define the penalized smoothed estimator as
\[
\widehat{\bbeta}_{\mathbb{S}}=\arg \min_{\mathbf b}\left\{\mathcal{R}_n(\bb) + \lambda_n \sum_{j=1}^p \sqrt{d} \rho \left( \| \bb_j^T \bR_j\|_{\gamma_j}  + \sqrt{\bb_j^T \bM_j\bb_j} \right)\right\}, \ \ \ \mbox{ for } \gamma_j \geq 2,
 \]
 for a convex and subadditive choice of  $\rho$. The smoothing matrix, $\bM_j \in \mathbb{R}^{d \times d}$, contains the inner products of the second derivatives
of the B-spline basis functions, i.e.,
$$\{\bM_j\}_{kl}=\int \Psi_k^{''}(x_j) \Psi_l^{''}(x_j) d x_j, \qquad  \bM_j=\bR_j^T \bR_j,$$
 $k,l=1,\cdots,d$, and  $\bR_j \in \mathbb{R}^{d \times d}$ is a matrix obtained from Cholesky decomposition of $\bM_j$.  
 Then, we can rewrite the problem as
\[
\widehat{\bbeta}_{\mbox{\sc \tiny s}}=\arg \min_{ \tilde{\mathbf b}}\left\{\mathcal{R}_n(\tilde{\bb})+ \lambda_n \sum_{j=1}^p  \sqrt{d}  \rho \left( \| \tilde{ \bb}_j \|_{\gamma_j}  +  \| \tilde{ \bb}_j \|_2 \right)\right\},
\]
with $\tilde {\bb}_j=\mathbf R_j \bb_j$ and
$$\mathcal{R}_n(\tilde{\bb})=- \frac{1}{n} \sum_{i=1}^n  \int_0^\tau \Bigl\{ \sum_{j=1}^p \tilde{\bb_j}^T \bR_j^{-1} \bPsi(X_{ij})+ \log \bigl( \frac{1}{n} \sum_{i=1}^n Y_i(t) \exp\{ \sum_{j=1}^p \tilde{\bb_j}^T \bR_j^{-1} \bPsi(X_{ij})\}\bigl)\Bigr\}d N_i(t). $$
  A crucial part of extending the previous results to this novel setting requires extending the results of Lemma \ref{lemma:min1} and Propositions \ref{lem:approx} and \ref{prop:approx} to the new penalty structure. Details of the proof are  presented in the Appendix E.

\begin{lemma} \label{lemma:smooth}
Equivalent  results to Lemma \ref{lemma:min1}  and  Propositions \ref{lem:approx}, \ref{prop:approx}   hold for $\mathcal{R}_n(\tilde{\bb})$,  
 on a event   that  has probability very close to 1 (details are presented in the Appendix). 
    \end{lemma}
With the help of the results presented in earlier sections and this Lemma, we have the following Corollary.

\begin{corollary}
Let conditions of Theorem 2 be satisfied. Let $\bM_j$ be well defined with

$\underline\lambda= \min_{1 \leq j \leq p}  \lambda_{\min}^{}(\bR_j) >0$.
Then, for some constant $A >4$ and the choice of the tuning parameters
$$ \lambda_n d^2 \geq  A  \min\left\{ \zeta^2, \sqrt{\frac{\log (pd)}{n}} \right\},$$
 with the probability of at least of $1- \delta$, $\delta>0$,
    \begin{equation}
\| f_{ \widehat{\boldsymbol\beta}_{ \mathbb S}} - f_{\boldsymbol\beta^*}\| _{ }^2
 \leq
{\underline {\underline \omega}}^{-1} e^{   C ( r_n e^{C}  + 2\upsilon_1  )} 32 \sqrt{2} \frac{s \lambda_n^2d }{\zeta^{2} } \sum_{j \in \mathcal{M}_*} \mathbf R_j \mathbf R_j^T,
\end{equation}
 for
  $r_n=   \zeta^{-2}s \lambda_{n} d$ , and $0 \leq \upsilon_1 \leq 1$  satisfying
 \begin{equation}\label{eq:upsilons2}
\upsilon_1 e^{-2 C \upsilon_1 }  \leq 16 \lambda_n^2\underline\lambda   \frac{  s d }{\zeta^2}  \rho'(0+) .
 \end{equation}

%
\end{corollary}

 To the best of our knowledge, this is the first finite sample result on prediction properties of a non-parametric smoothing estimator for the high-dimensional Cox model. A particular example of a smooth selection  is the Elastic net penalty \citep{ZH05}. Although our previous results easily apply to this penalty (by specifying $\gamma_j=1$ and $\rho=l_1$), its efficient implementation  in the Cox model was only recently proposed in \cite{W12}, but its theoretical properties have not been previously studied.
Although tackled as the last problem,  the importance of the obtained finite sample bounds for smooth selection lies in the inadmissibility of such results with techniques that already exist in the literature.  In particular, in the case of Elastic-Net penalty we obtain with high probability,
\[
\| f_{ \widehat{\boldsymbol\beta}_{ \mbox{elastic net}}} - f_{\boldsymbol\beta^*}\| _{ }^2
 \leq
{\underline {\underline \omega}}^{-1} e^{   C ( r_n e^{C}  + 2\upsilon_1  )} 32 \sqrt{2} \frac{s^2 \lambda_n^2  }{\zeta^{2} }  
\]
 for
  $r_n=   \zeta^{-2}s \lambda_{n} $ , and $0 \leq \upsilon_1 \leq 1$  satisfying $
\upsilon_1 e^{-2 C \upsilon_1 }  \leq 16 \lambda_n^2    \frac{  s   }{\zeta^2}  \rho'(0+) .
$

%


\section*{Discussion}

In this paper, we propose a new method for analyzing the theoretical oracle risk properties of likelihood functions that are not necessarily of a quadratic nature. By sandwiching  the likelihood with two other processes, we show that  it is sufficient to analyze the risk properties of the bounding processes alone.  
To the best of our knowledge, minimax rates,  have not been established for any survival model so far despite their importance.
Equivalents of traditional information theoretic tools, such are Fano's lemma, are not  easy to understand in the Cox model setup. Our proposed method of sandwiching the likelihood with two quadratic likelihoods may be useful in establishing minimax rates. 

%

\appendix


\section{The Restricted Eigenvalue Condition}\label{sec:RE}

The restricted eigenvalue condition, {\bf{RE}}$(\mu,s,\rho,\bgamma)$, defined in (12) represents a generalization of the cone constraint condition that appears in work on Lasso problems \citep{BRT09}.    Equivalent definitions were proposed for various hazard rate  models \citep{L12,GG11,KN12,HZ12}.
   We refer to \\ \cite{BvG11} for comparisons of different kinds of compatibility and restricted eigenvalue conditions and their relationships for sparse linear models.
  The usual scaling factor of $\sqrt{n}$ disappears in the definition of the restricted eigenvalue condition because it is included in the definition of the empirical norm, $\| f (\cdot)\|^2_{n,\cdot}$. Compared to the RE condition in \cite{BRT09}, the denominators differ in that the $l_2$ norm is replaced with an $l_{1,\boldsymbol\gamma}$ norm.
 In the least squares  procedures,  $ \bigtriangledown^2 \mathcal{L}_n(\bbeta^* )  = -\bX^T \bX$ and the restricted eigenvalue conditions are defined on the eigenvalues of $ \bX^T \bX$.  Condition \eqref{eq:RE} can be seen as a  rescaling of the minimum eigenvalue problem in the classical {\bf RE} condition needed for the complex likelihood structures.

Determining the class of matrices that satisfy the {\bf RE($\mu,s,\bgamma$)} condition is an important open question. Heuristically we can argue in the following manner. First, we observe that
with respect to time, $\int_0^\tau \bV_n(0,t) d \bar N(t)$   has a martingale  structure. With respect to    $\bbeta^*  $, it is  a  function of the  matrix  $\sum_{i=1}^n \sum_{q=1}^n\bPsi^T(\bX_i)\bPsi(\bX_q)$. Using Condition \ref{cond:survival} and the boundedness of the $\Psi$ functions, matrix $\int_0^\tau \bV_n(0,t) d \bar N(t)$ will belong to a random matrix ensemble with sub-gaussian tails, studied in \cite{Z09}. Dependence through time was shown not to be essential in  \cite{HZ12}, where a  lower bound for {\bf RE}  was shown to be independent of time.
Moreover, we can  combine  both  results to conclude  that for large enough
  sample size, there exists a positive constant  $\zeta_1$   such that with overwhelming probability
   \[
 \ \ \ \min_{{\boldsymbol\Delta} \in \mathbb{C}_{\mu,\rho}, {\boldsymbol\Delta} \neq 0  }  \frac{\| {\boldsymbol\Delta}^T\{-\bigtriangledown^2 \mathcal{L}_n(\boldsymbol\beta^*) \}^{}{\boldsymbol\Delta} \|_2}{ \| {\boldsymbol\Delta}_{\mathcal{M}_*} \|_{1,\boldsymbol\gamma}^2} \geq \zeta_2
   .\]

\section{Preliminary Lemmas}

The following lemma provided exponential inequality for a martingale sequence and  can be found in   \cite {vG95} as Lemma 2.1
\begin{lemma}\label{lem:expineq}
Let $(\Omega, \mathcal{F},{P})$ be a probability triple and let $M_t$ be a sequence of locally square integrable martingales w.r.t. the filtration $\mathcal{F}_t$. Suppose that $|M_t-M_{t-}|\leq K$ for all $t >0$ and some $0<K<\infty$. Then, for each $a>0,b>0.$
\[
P\left( M_t \geq a \mbox{ and } \langle M,M\rangle_t \leq b^2 {\mbox{ for some }t}\right) \leq \exp\left\{ - \frac{a^2}{2(aK +b^2)}\right\},
\]
where $ \langle M,M\rangle_t$ denotes predictable variation of the martingale sequence $M_t$.
\end{lemma}

The following lemma provides an exponential inequality for a unbounded  supermartingale sequence and can be found in \cite{FGL11} as Corollary 2.3.

\begin{lemma}\label{lem:nonnegineq}
Let $(\Omega, \mathcal{F},{P})$ be a probability triple . Assume that  $(\xi_i,\mathcal F_i)_{i=1,\dots,n}$ are supermartingale differences i.e. $E (\xi_i|\mathcal F_{i-1})\leq 0$. Let $b>0$ and 
\[
V_k^2(b)=\sum_{i=1}^k E \left( \xi^2_i 1\{\xi_i \leq b\} | \mathcal F_{i-1}\right), \qquad k=1,\dots, n.
\]
Then, for any $a \geq 0$, $b>0$ and $c>0$
\[
P \left( \sum_{i=1}^k \xi_i \geq a \mbox{ and } V_k^2(b) \leq c^2  {\mbox{ for some }k}\right) \leq  \exp\left\{ - \frac{a^2 }{2(c^2  + \frac{1}{3} ab)} \right\} + P(\max_{1 \leq i \leq n}\xi_i > b ).
\]
\end{lemma}

\section{Proofs of Propositions}

\begin{proof}[Proof of Proposition \ref{prop:w}]
Without loss of generality, let us represent the optimization problem \eqref{eq:opt} as a quadratically  constrained minimum of the ratio of two quadratic functions of the following form
\begin{equation}\label{eq:opt1}
\begin{array}{cc}
\min_{\mathbf b}  &\sum_{q=1}^N\frac{(\bb^T \bA_1^i \bb + 2 {\ba_1^i}^T \bb + c_1) 1(i \in \mathcal{R}_q)}{ \bb^T \bA_2 \bb + 2 \ba_2^T \bb + c_2},
\\
\ \\
\mbox{ s.t }& \| \bb \|_2^2 \leq r_n, \\
& \mathbf b \in \mathbb{R}^{pd}, \\
\end{array}
\end{equation}
where $\bA_1^i=\bPsi(\bX_i)\bPsi(\bX_i)^T  $, $\bA_2=\sum_{l \in \mathcal{R}_q} \bPsi(X_l) \bPsi^T(\bX_l)$ and $a_1^i= \bPsi(\bX_i)$, $a_2= \sum_{l \in \mathcal{R}_q} \bPsi(X_l)$. Constants $c_1$ and $c_2$ are residuals of the Maclaurin series expansions of the functions $\exp\{\bb^T \bPsi(\bX_i)\}$ and $\sum_{l \in \mathcal{R}_q} \exp\{\bb^T \bPsi(\bX_l)\} $. This makes $\bA_1^i$ and $\ba_1^i$ second order and first order approximations of $\exp\{\bb^T \bPsi(\bX_i)\}$, around $\mathbf 0$ .  

Condition \eqref{eq:optcond}  implies that for any feasible point $\bb$, the above optimization problem is well defined. Multiplying \eqref{eq:optcond} by $(\bb^T,1)$  from the left and $(\bb^T,1)^T$ from the right results in
\[
\sum_{l \in \mathcal{R}_q} \exp\{\bb^T \bPsi(\bX_l)\} + \eta(\|\bb\|_2^2 - r_n) \geq \delta(\| \bb\|_2^2)+1,
\]
which implies that $\sum_{l \in \mathcal{R}_q} \exp\{\bb^T \bPsi(\bX_l)\}  \geq \delta(\| \bb\|_2^2)+1 \geq \delta >0$.

Let us fix an $i \in \mathbb{R}_q$ for some $q$. Now, let us define
\begin{eqnarray}
d_1=\inf \bigl\{ f(\bb): \| \bb\|\leq r_n,  \bb^T \bA_1^i \bb + 2 {\ba_1^i}^T \bb + c_1 \geq 0\bigl\}, \\
d_2=\inf \bigl\{ f(\bb): \| \bb\|\leq r_n,  \bb^T \bA_1^i \bb + 2 {\ba_1^i}^T \bb + c_1 \leq 0\bigl\},
\end{eqnarray}
with  $f(\bb)={(\bb^T \bA_1^i \bb + 2 {\ba_1^i}^T \bb + c_1)}/{( \bb^T \bA_2 \bb + 2 \ba_2^T \bb + c_2)}$.
Then using the relation that
\[
\inf\{f(\bb): \bb \in \mathcal{C}_1 \cup \mathcal{C}_2\} = \min \bigl\{  \inf_{\mathbf b \in \mathcal{C}_1} f(\bb),  \inf_{\mathbf b \in \mathcal{C}_2} f(\bb) \bigl\},
\]
we have that the optimal solution to \eqref{eq:opt1} is equal to $\min\{d_1,d_2\}$. By definition, $d_1$ is nonnegative. It remains to show that $d_2$ is finite. Indeed, for every $\bb$ satisfying $\|\bb\|_2^2\leq r_n$ and $\bb^T \bA_1^i \bb + 2 {\ba_1^i}^T \bb + c_1\leq 0$, we have
\[
d_2 \geq f(\bb) \geq \frac{\bb^T \bA_1^i \bb + 2 {\ba_1^i}^T \bb + c_1}{\delta(\| \bb\|_2^2)+1} \geq \frac{1}{\delta}
\lambda_{\min} \left(\begin{array}{cc} \bA_1^i & \ba_1 \\ \ba_1^T & c_1\end{array} \right).
\]
\end{proof}

\begin{proof}[Proof of Proposition \ref{lem:approx}]
 To see that the equation \eqref{eq:sand1} is correct, we adopt the  following reasoning. First, note that $\| f_{\mathbf b} - f_{\boldsymbol\beta^*}\|_{n, \mathbf b^*}^2$ is equal to
\[
n^{-1}\int_0^\tau \frac{\sum_{i,q=1}^n w_i w_q ( a_i -  a_q)^{\otimes 2} e^{(1-c)a_i -\bar c} e^{(1-c)a_q -\bar c} }{\sum_{i,q=1}^n 2 w_i w_q e^{(1-c)a_i -\bar c} e^{(1-c)a_q -\bar c} } d \bar N(t),
\]
with $a_i=  (\bb - \bbeta^*)^T(\Psi(\bX_i) - \bE_n(\bbeta^*,t)) $ and $w_i=Y_i(t) \exp \{ {\bbeta^*}^T \Psi(\bX_i)\}$ and $\bar c= (1-c)(\max_i a_i + \min_i a_i)/2$. If we let $\eta = a_{ {\mathbf b} - \boldsymbol \beta^*}$, we can see that $ \max_i |(1-c) a_i -\bar c| \leq \eta/2$. ({bf{should the following be true for max over i?}}) Using this notation,  $
e^{(1-c) a_i -\bar c} \geq e^{- \eta/2} \mbox{ and } e^{(1-c) a_i -\bar c} \leq e^{ \eta/2}
$ leading to
\begin{eqnarray*}
\| f_{\mathbf b} - f_{\boldsymbol\beta^*}\|_{n, \mathbf b^*}^2 &\geq& \exp \{ - 2 \eta \} n^{-1} \int_0^\tau \frac{\sum_{i,q=1}^n w_i w_q ( a_i -  a_q)^{\otimes 2}  }{\sum_{i,q=1}^n 2w_i w_q  } d \bar N(t) \\
&=& \exp \{ - 2 \eta \} \ \| f_{\mathbf b} - f_{\boldsymbol\beta^*}\|_{n, \boldsymbol\beta^*}^2.
\end{eqnarray*}
The upper bound follows the same reasoning, and thus it is omitted.
The lower bound of the RHS of previous inequality  follows by repeating the same steps as in  Proposition \ref{prop:approx} and the definition of the weight vectors, $\omega_i(\bbeta^*)$, in \eqref{eq:weights},
  \[
  \| f_{\mathbf b} - f_{\boldsymbol\beta^*}\|_{n, \mathbf b^*}^2 \geq  \underline{\underline \omega} \exp \{ - 2 \eta \} \ \| f_{\mathbf b} - f_{\boldsymbol\beta^*}\|_{n}^2.
  \]
The upper bound follows directly from Proposition \ref{prop:approx}  by taking $\bb^*=\bbeta^*$ to obtain 
\[
\| f_{\mathbf b} - f_{\boldsymbol\beta^*}\|_{n, \mathbf b^*}^2  \leq  \exp \{ - 2 \eta \} \ \| f_{\mathbf b} - f_{\boldsymbol\beta^*}\|_{n}^2.
\]
%

\end{proof}

\begin{proof}[Proof of Proposition \ref{prop:approx}]
Let $N$ denote  the cardinality of the set $\{i=1,\cdots,n: N_i(\tau) =1\}$.
The  weight process, $\omega_i (\bb,t)$ as defined in \eqref{eq:1}, satisfies the following normalization  uniformly over $\bb$ and $t$,
$$
\frac{1}{n}\sum_{i=1}^n Y_i(t) \omega_i (\bb,t) =1.
$$
For each $\bb$, there exists at least one $i\in\{1,\dots,n\}$ such that $\omega_i(\bb,t) >0$ and that for all $i$,  for which $\exists t \in[0,\tau]$, $Y_i(t)=1$, we have that $\omega_i(\bb,t)\leq n$, for all $t$.  

 %
Let us  denote
$$
\omega_i(\bb)=\int_0^\tau Y_i(t) \omega_i(\bb,t) d \bar N(t),
$$
with $\omega_i(\bb,t)$ defined as in \eqref{eq:1}. 
If $t_1<\dots<t_N$ are ordered failure times and $\mathcal{R}_j=\{i\in\{1,\dots,n\}: Z_i \geq t_j\}$ is at risk set, then $\omega_i(\bb)$ has the following representation:  
$$
\omega_i(\bb)
=  \sum_{j=1}^N  \frac{\exp\{\bb^T \bPsi(\bX_i)\} 1\{i \in \mathbb{R}_j\}}{ \sum_{l \in \mathbb{R}_j} \exp\{\bb^T \bPsi(\bX_l)\}},
$$
which matches the definition provided in Theorem \ref{cor:localSOI}  equation \eqref{eq:weights}.
Note that $\omega_i \geq 0$ and $\omega_i >0$ for $  i \in \{1,\dots,n\}   $.
 Using the previous notations, we have
\[
\| f _{\mathbf b}\| _{n, \mathbf b^*}^2 = \frac{1}{n} \sum_{i=1}^n f_{\mathbf b}^2(X_i) \omega_i(\mathbf b^*) -  \left(\frac{1}{n} \sum_{i=1}^n f_{\mathbf b}(X_i) \omega_i(\mathbf b^*)\right)^2,
\]
With this notation at hand, we have that
\[
\frac{1}{n} \sum_{i=1}^n\omega_i(\mathbf b^*)=\frac{1}{n} \sum_{j=1}^N   \sum_{i \in \mathbb{R}_j} \frac{\exp\{\bb^T \bPsi(\bX_i)\}}{ \sum_{l \in \mathbb{R}_j} \exp\{\bb^T \bPsi(\bX_l)\}} = \frac{N}{n}.
\]
Since  $\omega_i(\bb^*)\geq 0$,  and are defined as conditional probabilities  we have $ \omega = \max\{\omega_i(\bb):   i \in \{1,\dots,n\} ,  \mathbf b \in R^{pd}\}   \leq 1$.
We are then able to conclude that $1 \geq \bar  \omega = \max\{\omega_i(\bb):   i \in \{1,\dots,n\} ,  \mathbf b \in R^{pd}\} \geq 1/ n \geq \underline  \omega = \min\{\omega_i(\bb):  i \in \{1,\dots,n\} ,  \mathbf b \in R^{pd}\}$  .
Hence,
\[
\| f _{\mathbf b}\| _{n, \mathbf b^*}^2\leq  \bar \omega \| f_{\mathbf b} \|_n^2 - \underline \omega \left(\frac{1}{n} \sum_{i \in \mathbf I} f_{\mathbf b}(X_i) \right)^2 \leq   \| f_{\mathbf b} \|_n^2.
\]
To obtain the left-hand side of \eqref{eq:approx},  remember from previous exposition we have  \begin{eqnarray*}
\| f _{\mathbf b}\| _{n, \mathbf b^*}^2
&=& \frac{1}{n}\sum_{i=1}^n  \omega_i ({\mathbf b}_{\mathbf b})   (f_{\mathbf b}(\bX_i)-\bar f_{\mathbf b}^{*})^2
\end{eqnarray*}
with
$
\bar f_{\mathbf b}^{*} = \frac{1}{n}\sum_{i=1}^n \omega_i ({\mathbf b}_{\mathbf b})  f_{\mathbf b}(X_i)
$
and  $\omega_i ({\mathbf b}_{\mathbf b}) $ following the definition in \eqref{eq:weights}.
Hence, by centering  the data so that the sample mean is  equal to zero, i.e., $\frac{1}{n}\sum_{i =1}^n  f_{\mathbf b}(\bX_i)=0$, we have  
\begin{eqnarray*}
\| f _{\mathbf b}\| _{n, \mathbf b^*}^2 &\geq&\underline \omega \frac{1}{n}\sum_{i \in I}  \left( f_{\mathbf b}^2(\bX_i)+\{\bar f_{\mathbf b}^{*}\}^2 \right)  +  2 \underline \omega \ \bar f_{\mathbf b}^{*} \  \left(\frac{1}{n}\sum_{i =1}^n  f_{\mathbf b}(\bX_i)\right)\\
&\geq &\underline \omega \frac{1}{n}\sum_{i =1}^n   f_{\mathbf b}^2(\bX_i) = \underline \omega \|   f_{\mathbf b}\|_n^2.
\end{eqnarray*}

The result of the   Proposition \ref{prop:approx}
 follows easily after applying Proposition \ref{prop:w} on the interval $[-b_n,b_n]$ and following discussion after Proposition \ref{prop:w}.
 \end{proof}

  \section{Proofs of Lemmas}

  \begin{proof}[Proof of Lemma \ref{lemma:min1}]
  Let us first concentrate on the first statement   of  \eqref{eq:temp11}.
This can be seen from the following reasoning. Let us define a function
$$
f(\bb):= -(\bb - \bbeta^*)^T\bv_n  +  \lambda_n P(\bb) - \mathcal{L}_n(\bbeta^*),
$$
where $\bv_n \in \mathbb{R}^{pd}$.
First, we show that zero is a local minimum of  function $f(\bb)$ for all $\bb$ such that $\|\bb_j\|_{1} \leq 1$.
Note that
$$
f(\bb) -f(\mathbf 0) = \sum_{j=1}^p \biggl( -\bb_j^T \bv_{n,j} +\lambda_n d^{1/{\gamma_j^*}}\rho(\| \bb_j\|_{\gamma_j}) \biggl),
$$
 and conditional on the event $\mathcal{E}_{n,j}=\left\{ \|\bv_{n,j} \|_{\gamma_j^*} \leq \lambda_n d^{1/{\gamma_j^*}} {\rho'(0+)}  \right\}$,
\begin{eqnarray*}
-\bb_j^T \bv_{n,j}  +\lambda_n d^{1/{\gamma_j^*}}\rho(\| \bb_j\|_{\gamma_j}) \geq \| \bb_j\|_{\gamma_j} \left(  -\|\bv_{n,j} \|_{\gamma_j^*} + \lambda_n d^{1/{\gamma_j^*}}{\rho'(0+)} \right)\geq 0,
\end{eqnarray*}
where we have utilized the H\"{o}elder inequality. Therefore, we can conclude that $f(\bb) -f(\mathbf 0)  \geq 0$ if the event $\mathcal{E}_n = \cap_{j=1}^p \mathcal{E}_{n,j}$. Because $f$ is a convex function, we can conclude that $0$ is a global minimum as well. Note that we don't require unicity of minimum.

We are left to prove the second statement of \eqref{eq:temp11}.
We proceed in the similar way by first defining an appropriate function to minimize over. To that end, let us define
$$
f(\bb):= -|( \bbeta^*-\bb )^T\bv_n | +  \lambda_n P(\bb) - \mathcal{L}_n(\bbeta^*  ),
$$
where $\bv_n \in \mathbb{R}^{pd}$.
By the same reasoning as above, it suffices to notice that
\begin{eqnarray*}
&&- \left( |(\bbeta^*-\bb)^T \bv_{n,j}|- |{\bbeta^*}^T \bv_{n,j}|  \right)+\lambda_n d^{1/{\gamma_j^*}}\rho(\| \bb_j\|_{\gamma_j})  \\
&\geq&  - | \mathbf b^T \bv_{n,j}|  +\lambda_n d^{1/{\gamma_j^*}}\rho(\| \bb_j\|_{\gamma_j})  \\
&\geq& \| \bb_j\|_{\gamma_j} \left(  -\|\bv_{n,j} \|_{\gamma_j^*} + \lambda_n d^{1/{\gamma_j^*}}{\rho'(0+)} \right)\geq 0,
\end{eqnarray*}
by first using $|x-y| \geq | |x| - |y||$
 and then H\"{o}elder inequality.
  \end{proof}

%
%

  \begin{proof}[Proof of Lemma \ref{lemma:temp1}]

We make use of the following decomposition 
$$\| \bE_n(\bbeta^*,t)-\mathbf e(\bbeta^*,t) \|_\infty  = \hskip 270pt$$
\begin{equation} \label{eq:temp100}
\begin{array}{ccc}
&& \leq  \max_{1\leq j \leq p, 1 \leq k \leq d} \frac{ \left|  {\{S_n^{(1)}\}_{jk}(\bbeta^*,t)  - \{s^{(1)}\}_{jk}(\bbeta^*,t) }\right|}{|s^{(0)}(\bbeta^*,t)| } \\
&+&  \max_{1\leq j \leq p, 1 \leq k \leq d} \left| \{s^{(1)}\}_{jk}(\bbeta^*,t)\right|  \left| \frac{1 }{S_n^{(0)}(\bbeta^*,t) } -\frac{1 }{s^{(0)}(\bbeta^*,t) }\right|
\end{array} := I_1 + I_2
\end{equation}
We will prove maximal inequalities for each of the two terms in the above inequality.

First,  consider  classes of functions indexed by $t$:


$$\mathcal{F} = \{1\{z>t\} \exp\{f_{\boldsymbol\beta^*}(x)\} / u  : t\in[0,\tau]\},
$$
and
$$\mathcal{G}^k = \{1\{z>t\}  \Psi_k(x) \exp\{f_{\boldsymbol\beta^*}(x)\} / u  : t\in[0,\tau] \}. $$

Since $\bbeta^*$ is a $s$-sparse vector we have that
$u=\exp\{\sum_{j \in \mathcal{M}_*} \| \bbeta_j^* \|_1\}$. We proceed by calculating theirs bracketing number. Noticing that previous classes are products of a class of indicator functions and a class of bounded functions  we have that
\[
\mathcal{N}_{[]}({\epsilon,\mathcal{F},L_2} )\leq 2/\epsilon^2, \qquad \mathcal{N}_{[]}({\epsilon,\mathcal{G}^k,L_2} )\leq 2/\epsilon^2,
\]
By direct consequence of theorem 2.14.9 of Van der Vaart and Wellner (1996) we obtain that there exists a constant $W$ such that
\[
P\left( \sqrt{n} \sup_{t \in [0,\tau] } \left| \frac{1}{n} \sum_{i=1}^n Y_i(t) \exp\{f_{\boldsymbol\beta^*}(\bX_i)\} / u -  E_{Y,X} Y_i(t) \exp\{f_{\boldsymbol\beta^*}(\bX_i)\} / u\}\right| \geq r\right)   \leq \frac{1}{2e} W^2 e^{-r^2}
\]
and
\begin{eqnarray*}
P\biggl( \sqrt{n} \sup_{t \in [0,\tau] } \biggl| \frac{1}{n} \sum_{i=1}^n Y_i(t) \Psi_k(\bX_i)\exp\{f_{\boldsymbol\beta^*}(\bX_i)\} / u   \biggl. \biggl. \hskip 80pt
\\
\hskip 50pt -  \biggl. \biggl. E_{Y,X} Y_i(t)\Psi_k(\bX_i) \exp\{f_{\boldsymbol\beta^*}(\bX_i)\} / u\}\biggl| \geq r\biggl)   \leq \frac{1}{2e} W^2 e^{-r^2},
\end{eqnarray*}
for every fixed $k \in \{1,\cdots, d\}$.
By replacing $r$ with $\sqrt{n} r_n$ in the first and utilizing union bound and replacing $r$ with $\sqrt{n r_n^2 + \log 2d}$ in the second we obtain
\begin{equation}\label{eq:proof1}
P\left( \sup _{ t \in [0,\tau]} \left | S_n^{(0)}(\bbeta^*,t) - s^{(0)}(\bbeta^*,t) \right| \geq u r_n\right) \leq \frac{1}{2e} W^2  e^{- nr_n^2},
\end{equation}
 
 \begin{equation}\label{eq:proof2}
 P\left( \sup _{ t \in [0,\tau]}  \| S_n^{(1)}(\bbeta^*,t) - s^{(1)}(\bbeta^*,t) \| _\infty \geq u \left( \sqrt{r_n^2 + \frac{\log 2d}{n}} \right)\right) \leq \frac{1}{4de}  W^2  e^{-n r_n^2},
\end{equation}

Second, from the definition of $s^{(0)}(\bbeta^*,t)$ and Condition \ref{cond:survival} (iii) we  observe that there exists a constant $0<D <1$  with  $D = P( Y(\tau) =1 )$ 
and
\[
\inf_{t \in[0,\tau]}   \frac{1}{n} \sum_{i=1}^n E_{Y,X} Y_i(t) \exp\{f_{\boldsymbol\beta^*}(\bX_i)\}  \geq \exp\{- m^* C\} P( Y(t) =1 ) >D\exp\{- m^* C\}
\]
with $C$ being an upper bound on $|\Psi_k(x)|$ and $m^*$ defined as minimum signal strength in  the additive component of the hazards model \eqref{eq:hazard}.

According to \eqref{eq:temp100}  and \eqref{eq:proof1} we have
\[
I_2 \leq  \frac{\sup_{t \in [0,\tau]} \left\|  s^{(1)} (\bbeta^*,t)\right\|_\infty  } {D\exp\{- m^* C\}} \frac{u r_n }{ \inf_{t \in [0,\tau]}S_n^{(0)}(\bbeta^*,t) }
\]
with probability  $\frac{1}{2ed} W^2  e^{- r_n^2}$, and according to  \eqref{eq:temp100}  and \eqref{eq:proof2} we have
\[
I_1 \leq \frac{u \left( \sqrt{ r_n^2 +\frac{\log 2d }{n}} \right) \exp\{m^* C\}}{D} \leq \frac{u \left(r_n + \sqrt{\frac{\log 2d }{n}} \right) \exp\{m^* C\}}{D} ,
\]
with probability  no smaller than $1-\frac{1}{4e}  W^2  e^{-n r_n^2}$.
To further bound $I_2$
  we show that  $  |S_n^{(0)} (\bbeta^*,t)|  $  is bounded away from zero with high probability.
To that end, we  employ Massart's Dvoretzky-Kiefer-Wolfowitz  inequality  bounding   how close an empirically determined distribution function is to the distribution function from which the empirical samples are drawn. Hence, {???}
\begin{eqnarray}\nonumber
&&P \left( \frac{1}{n} \sum_{i=1}^n \mathbbm{1}(Z_i \geq \tau) \geq \frac{1}{2}P(Z_1 \geq \tau) \right)
\\\nonumber
&&\geq  P \left( \sup_{t \in [0,\tau]} \sqrt{n} \left| \frac{1}{n} \sum_{i=1}^n \mathbbm{1}(Z_i \geq t)  - P(Z_1 \geq \tau) \right| \leq \sqrt{n} /2\  P(Z_1 \geq \tau) \right)
\\
&&\geq 1- 2 e^{-  n D^2 /2 }.
\label{eq:dvoretsky}
\end{eqnarray}
 
Remeber that $S_n^{(0)} (\bbeta^*,t) = \frac{1}{n} \sum_{i=1}^n \mathbbm{1}\{Z_i \geq t\} \exp\{f_{\boldsymbol \beta^*}(\bX_i)\}$
 and observe that for all $t\leq \tau$ we have $\{Z_i \geq t\} \supset\{ Z_i \geq \tau\}$.  Hence, 
\[
S_n^{(0)} (\bbeta^*,t) \geq \exp\{-m^* C\} \frac{1}{n} \sum_{i=1}^n \mathbbm{1}\{Z_i \geq \tau\} , \qquad \mbox{for all } t \leq \tau.
\]
Together with \eqref{eq:dvoretsky} we have
\[
P\left(\inf_{t \in [0,\tau]}S_n^{(0)} (\bbeta^*,t)  \geq \exp\{-m^* C\} D/2 \right) \geq \left( \frac{1}{n} \sum_{i=1}^n \mathbbm{1}(Z_i \geq \tau) \geq  D/2 \right) \geq 1- 2 e^{-  n D^2 /2 }.
\]
Next, we bound $\sup_{t \in [0,\tau]} \| s^{(1)}(\bbeta^*,t)\|_\infty $. Observe that 
$$\sup_{t \in [0,\tau]} \| s^{(1)}(\bbeta^*,t)\|_\infty \leq \sup_{t \in [0,\tau]} E_{X}( P\{Z_1 \geq t | \bX_1\} 
\exp\{ f(\boldsymbol \beta^* (\bX_1))\})
$$
$$\leq E_{X}(  \exp\{{\boldsymbol \beta^*} ^T \bPsi(\bX_1)\}) \leq  \exp\{C \log u \}  
 $$

With all of the above notice that  
\[
 \qquad I_2 \leq     \frac{2 u   r_n  \exp\{2m^* C\}  \exp\{C \log u \}  }{D^2}
\]
  with probability no smaller than $1-\frac{1}{2ed}  W^2  e^{-n r_n^2}  - 2 e^{-  n D^2 /2 }$ .
Hence, we conclude that 
\[
\| \bE_n(\bbeta^*,t)-\mathbf e(\bbeta^*,t) \|_\infty  \leq  \frac{u \left(r_n + \sqrt{\frac{\log 2d }{n}} \right) \exp\{m^* C\}}{D} +\frac{2 u   r_n  \exp\{2m^* C\}  \exp\{C \log u \}  }{D^2},
\]
with probability no smaller than $1-\frac{3}{8ed}  W^2  e^{-n r_n^2}  -  e^{-  n D^2 /2 }$.

\end{proof}

\begin{proof}[Proof of Lemma \ref{lemma:eventbounds}]
{\bf   Bounding $\mathcal{D}_{n,i}^c$}
Recall that
$$
\mathcal{D}_{n,i}= \bigcap_{j=1}^p\left\{ 4 \lambda_0(\tau)\left|    \int_0^\tau S_n^{(0)}(g,t)  dt \right|    \|  \bPsi(X_{ij})   \|_{\gamma_j^*}\leq \lambda_n  d^{1/\gamma_j^*}\rho'(0+)  \right\}.
$$
By simple union bound we see that 
\begin{equation}\label{eq:dn1}
P(\mathcal{D}_{n,i}^c) \leq \sum_{j=1}^p P \left( \lambda_0(\tau)\left|    \int_0^\tau S_n^{(0)}(g,t)  dt \right|    \|  \bPsi(X_{ij})   \|_{\gamma_j^*}\geq \lambda_n  d^{1/\gamma_j^*}\rho'(0+)  \right) .
\end{equation}



First,  observe that the definition of $S_n^{(0)} (g,t) $ allows the following bound
\[
\biggl|\int_{0}^{\tau}  S_n^{(0)} (g,t) d \Lambda_0(t) \biggl|\leq \Lambda_0(\tau)  \frac{1}{n} \sum_{i=1}^n \exp\{g(\bX_{i})\}  \int_{0}^{\tau}  Y_i(t) dt \leq \tau \Lambda_0(\tau)  \frac{1}{n} \sum_{i=1}^n \exp\{g(\bX_{i})\} ,
\]
whereas the boundedness of $\Psi_k$ allows
$
 \|  \bPsi(X_{ij})   \|_{\gamma_j^*} = \left( \sum_{k=1}^d \Psi_k^{\gamma_j^*}(X_{ij})\right)^{1/\gamma_j^*} \leq d^{1/\gamma_j^*} C 
$
to hold.
With this in mind,  we observe that
 \begin{equation}\label{eq:dn2}
P(\mathcal{D}_{n,i}^c) \leq \sum_{j=1}^p P \left( \tau\Lambda_0(\tau)  \lambda_0(\tau)    C^{1/\gamma_j^*}  \frac{1}{n} \sum_{i=1}^n \exp\{g(\bX_{i})\}  \geq \lambda_n   \rho'(0+)  \right) .
\end{equation}
Previous inequality is a tail probability of a sum of i.i.d.positive random variables where $g$ is the unknown function of interest.
By large-deviation inequality of  non-negative  random variables  (Lemma \ref{lem:nonnegineq} in the Appendix B), we obtain 
\begin{eqnarray}\label{eq:LargeDev}
P \left(  \sum_{i=1}^n \exp\{g(\bX_{i})\}
  \geq    \sqrt{n} \gamma_n  \right) \leq  e^{ - \frac{n\gamma_n^2}{2 \theta^2 + 2\gamma_n y/3}} + P(\max_{1 \leq i \leq n} \exp\{g(\bX_i)\} >y) ,
\end{eqnarray}
for a  sequence of non-negative numbers $\gamma_n$ and  a  truncation value $y$ such that 
\begin{eqnarray}\label{eq:theta}
\theta^2 \geq \sum_{i=1}^n E{ \exp\{2g(\bX_i) \}1\{ \exp\{2g(\bX_i) \} \leq y\}}.
  \end{eqnarray}
By choosing $\gamma_n = M \sqrt{n} \lambda_n \rho'(0+)$ with $M = 1/(\tau \lambda_0(\tau) \Lambda_0(\tau ) C)$, we obtain that 
\[
P(\mathcal{D}_{n,i}^c ) \leq e^{ - \frac{n^2 M^2 \lambda_n^2\rho'(0+)^2}{2 \theta^2 + 2 M\sqrt{n} \lambda_n \rho'(0+)y/3}} + P(\max_{1 \leq i \leq n} \exp\{g(\bX_i)\} >y).
\]

{\bf   Bounding $ \mathcal{E}_{n}^c$}

Notice that the set of interest, $ \mathcal{E}_{n}^c$, is a subset of
\[
 \bigcup_{j=1}^p \left\{ \|\bh_{n,{j}}(\bbeta^*)\|_{\infty} \geq  \lambda_n d^{1/\gamma_j^*} \rho'(0+)  \right \},
\]
 where $\|\bh_{n,{j}}(\bbeta^*)\|_{\infty} = \max_{1\leq k \leq d} | \{\bh_n\}_{jk}(\bbeta^*)|$.
 According to the definition  $\mathbf{h}_n(\bbeta^*) $ 
\begin{eqnarray}\label{eq:risk40}
\bh_n( \bbeta^*) = - n^{-1}\sum_{i=1}^n \int_0 ^\tau \left( \bE_n( \bbeta^*,t) -   \bPsi (\bX_i) \right) dM_i(t),
 \end{eqnarray}
   we   decompose $\mathbf{h}_n(\bbeta^*)  $  as follows
\begin{equation}\label{eq:51}
\mathbf{h}_n(\bbeta^*) :=\upsilon +\nu
\end{equation}
$$ = \frac{1}{n} \sum_{i=1}^n \int_0^\tau \left(\mathbf{E}_n(\bbeta^*)  - \be   (\bbeta^*,t) \right)   dM_i(t) + \frac{1}{n} \sum_{i=1}^n \int_0^\tau \left(\be (\bbeta^*,t) - \bPsi(\bX_i) \right)   dM_i(t)  .$$
 We will consider each term separately.
First, to control $\upsilon_{jk}$'s we develop a finite sample result  in Lemma \ref{lemma:temp1} whose proof   can be found  in  the Appendix E .

Next, we bound  $| {\boldsymbol\Delta} \upsilon_{jk} |$ and the predictable variation of the martingale $\upsilon_{jk}$. By Lemma \ref{lemma:temp1}, with high probability, the jumps are bounded by
\begin{eqnarray}\label{eq:jumps1}
 \hskip 30pt | {\boldsymbol\Delta} \upsilon_{jk} | = \frac{1}{n}  \left|\{\bE_n(\bbeta^*)\}_{jk} -\{\be(\bbeta^*)\}_{jk} \right|
\leq \frac{1}{n}  \sup_{0 \leq t \leq \tau} \left\|\{\bE_n(\bbeta^*,t)\} -\{\be(\bbeta^*,t)\} \right \|_\infty  \leq \frac{w_n}{n},
\end{eqnarray}
with $w_n=cr_n + \sqrt{\frac{\log d}{nu^2}}$.
The predictable variation process can be bounded as follows
\begin{eqnarray*}
\langle {\boldsymbol\Delta} \upsilon_{jk} \rangle_2  &=& \frac{1}{n^2}  \int_0^\tau \left[  \{\bE_n(\bbeta^*,t)\}_{jk} -\{\be(\bbeta^*,t)\}_{jk} \right]^2 d \langle\bar M(t)\rangle
\\
&\leq&\frac{1}{n}  \sup_{0 \leq t \leq \tau} \left\|\{\bE_n(\bbeta^*,t)\} -\{\be(\bbeta^*,t)\} \right \|_\infty^2 \int_0^\tau S_n^{(0)}(g,t) d \Lambda_0(t).
\end{eqnarray*}
The first term on the RHS of the above equation can be bounded above with high probability using Lemma \ref{lemma:temp1} with $w_n$. For the last term we  
use the result in \eqref{eq:LargeDev} to conclude that 
\begin{eqnarray}\label{eq:jumps2}
\langle {\boldsymbol\Delta} \upsilon_{jk} \rangle_2  
 &\leq& \frac{ \tau \Lambda_0(\tau)  }{n \sqrt{n}} w_n^2    \gamma_n  ,
\end{eqnarray}
for a sequence of non-negative numbers $\gamma_n$, with probability  larger than or equal to
\[
1- e^{ - \frac{n\gamma_n^2}{2 \theta^2 + 2\gamma_n y/3}} - P(\max_{1 \leq i \leq n} \exp\{g(\bX_i)\} >y)
\]
for any truncation value $y$ satisfying \eqref{eq:theta}.

Then, observe that  for any three events $A_1,A_2,A_3$, 
$$P(A_1)=P(A_1\cap A_2) + P(A_1|A_2^c) P(A_2^c) \leq P(A_1 \cap A_2) + P(A_2^c)$$ and  similarly 
$P(A_1 \cap A_2) \leq P(A_1 \cap A_2 \cap A_3) + P(A_3^c), $
leading to
$$
P(A_1) \leq P(A_1 \cap A_2 \cap A_3)+P(A_2^c)+P(A_3^c).
$$
Let $A_1=\{ |\upsilon_{jk}| \geq q_n \}$,$A_2 = \{| {\boldsymbol\Delta} \upsilon_{jk} |  \leq \frac{w_n}{n}\}$ and $A_3=\{\langle {\boldsymbol\Delta} \upsilon_{jk} \rangle_2 \leq \frac{ \tau \Lambda_0(\tau)  }{n \sqrt{n}} w_n^2    \gamma_n\}$.
 By   large deviation inequality for martingales of bounded jumps and variation in Lemma \ref{lem:expineq}, there exists a sequence of positive numbers $q_n$ such that
\[
P \left( |\upsilon_{jk}| \geq q_n  \right) \leq 2 e^{ -  \frac{n q_n^2 }{   K q_n + K_1^2  } } + P\Bigl(| {\boldsymbol\Delta} \upsilon_{jk} |  \geq \frac{w_n}{n}\Bigl) + P\Bigl(\langle {\boldsymbol\Delta} \upsilon_{jk} \rangle_2 \geq \frac{ \tau \Lambda_0(\tau)  }{n \sqrt{n}} w_n^2    \gamma_n\Bigl).
\] 
 By Lemma \ref{lemma:temp1} and equations \eqref{eq:jumps1} and \eqref{eq:jumps2} we have
\[
P \left( |\upsilon_{jk}| \geq q_n  \right) \leq 2 e^{ -  \frac{n q_n^2 }{   K q_n + K_1^2  } }  + \frac{3W^2}{8ed}    e^{-\frac{n r_n^2D^2}{u^2 e^{2m^*C}}}  + e^{- \frac{ n D^2}{2} }  + e^{ - \frac{n\gamma_n^2}{2 \theta^2 + 2\gamma_n y/3}} +P(\max_{1 \leq i \leq n} \exp\{g(\bX_i)\} >y) ).
\]
for $K=w_n/n$ and $K_1^2 = \gamma_n w_n^2 \tau \Lambda_0(\tau)/n\sqrt{n}$.
The choice of $\gamma_n$ is driven by \eqref{eq:LargeDev} where we considered  
$\gamma_n = M \sqrt{n} \lambda_n \rho'(0+)$ with $M = 1/(\tau \lambda_0(\tau) \Lambda_0(\tau ) C)$.
For a $q_n =\frac{1}{2} \lambda_n d^{1/\gamma_j^*} \rho'(0+) $, $K q_n \leq K_1^2$ as long as 
 \[
2     \omega_n   \geq C \lambda_0(\tau)  d^{1/\gamma_j^*} .
\]
With $w_n=cr_n + \sqrt{\frac{\log d}{nu^2}}$
  the choice of $r_n= C \lambda_0(\tau) \sqrt{n} d^{1/\gamma_j^*}   \sqrt{\frac{\log d}{ u^2}}$, suffices to guarantee the above inequality.
 For such choices of $\omega_n,\gamma_n$ and $q_n$ we have
\begin{eqnarray}\label{eq:temp17a}
P \left( |\upsilon_{jk}| \geq \frac{1}{2} \lambda_n d^{1/\gamma_j^*} \rho'(0+)   \right) \leq  \hskip 150pt \\
2 e^{ -  \frac{n q_n^2 }{   2 K_1^2  } } +\frac{3W^2}{8ed}    e^{-\frac{n r_n^2D^2}{u^2 e^{2m^*C}}}  +  e^{-   \frac{nD^2 }{2} } +e^{ - \frac{n\gamma_n^2}{2 \theta^2 + 2\gamma_n y/3}} +P(\max_{1 \leq i \leq n} \exp\{g(\bX_i)\} >y) ). \nonumber
\end{eqnarray}
The right-hand-side of (\ref{eq:temp17a}) can be simplified to 
\[
e^{-\frac{n^2 C \lambda_n \rho'(0+)}{2 \lambda_0(\tau)}}  + \frac{3W^2}{8ed}    e^{- \frac{n^2 C \lambda_0(\tau) D^2 d^{2/\gamma_j^*} \log d}{u^4 e^{2 m^* C}}}  +  e^{-  \frac{nD^2 }{2} } +e^{ - \frac{n^2 M^2 \lambda_n^2\rho'(0+)^2}{2 \theta^2 + 2 M\sqrt{n} \lambda_n \rho'(0+)y/3}} +P(\max_{1 \leq i \leq n} \exp\{g(\bX_i)\} >y)),
\]
which can be further bounded by
\[\leq \left(3+ \frac{3W^2 }{8ed}\right)e^{-n^2 {C}_{\lambda_n,n,p,d}} +P(\max_{1 \leq i \leq n} \exp\{g(\bX_i)\} >y)), 
\]
for 
\[
{C}_{\lambda_n,n,p,d}=\min \left\{ \frac{ C \lambda_n \rho'(0+)}{2 \lambda_0(\tau)},  \frac{  C \lambda_0(\tau) D^2 d^{2/\gamma_j^*} \log d}{u^4 e^{2 m^* C}}, \frac{D^2 }{2n},  \frac{ M^2 \lambda_n^2\rho'(0+)^2}{2 \theta^2 + 2 M\sqrt{n} \lambda_n \rho'(0+)y/3} \right\}
\]

Second, to control the $\nu$ term in \eqref{eq:51}, we observe that according to Lemma \ref{lemma:temp1}, there exists a constant $0<D = P(Y(\tau)=1) \leq 1$   such that for the $u$ as defined in Condition \ref{cond:survival} (iii) we have
\[
\sup_{ t \in [0,\tau]} \| \mathbf e (\bbeta^*,t)\|_\infty \leq \frac{C \sup_{ t \in [0,\tau]}  s^{(0)}(\bbeta^*,t)  }{ D \exp\{-m^* C\}  } \leq C u.
\]
Thus, each  $\nu_{jk}/u$ is a sum of a sequence of i.i.d bounded random variables. However, across $k$'s, i.e., group elements, $\nu_{jk}/u$ are not independent random variables.  By Hoeffding's inequality,
\begin{equation*}\label{eq:temptemp17}
P \left( \max_{1 \leq k \leq d}| \nu_{jk} | \geq   2 \| M\|_n C  u  t_n\right) \leq 2 e^{ -n t_n^2},
\end{equation*}
where  $\|M\|_n$ is proportional to $ E \sqrt{ \frac{1}{n} \sum_{i=1}^n M_i^2(\tau)}$. Because $\bar M$ is a bounded martingale, we  can conclude that there exists a constant $c_1 >0$ such that $ \| M\|_n \leq c_1$.

Hence, for $t_n = \lambda_n d^{1/\gamma_j^*} \rho'(0+) /4c_1 C u$ we obtain
\begin{equation}\label{eq:temptemp18}
P \left( \max_{1 \leq k \leq d}| \nu_{jk} | \geq  \frac{1}{2} \lambda_n d^{1/\gamma_j^*} \rho'(0+)  \right) \leq 2  e^{ -n \frac{ \lambda_n^2 d^{2/\gamma_j^*} {\rho'}^2(0+)}{16 c_1^2 C^2 u^2} }.
\end{equation}

Utilizing \eqref{eq:temp17a} and   \eqref{eq:temptemp18}   we obtain a   bound on the size of the set $ \mathcal{E}_{n}^c$ as follows,
$P \left( \mathcal{E}_{n}^c \right)$
\begin{eqnarray*}
&\leq&   2pd  \left (\max\biggl\{   \left(3+ \frac{3W^2 }{8ed}\right)e^{-n^2 {C}_{\lambda_n,n,p,d}} ,    e^{   -n \frac{ \lambda_n^2 d^{2/\gamma_j^*} {\rho'}^2(0+)}{16 c_1^2 C^2  u^2} } \biggl\}+P(\max_{1 \leq i \leq n} \exp\{g(\bX_i)\} >y)  \right)
 .  
\end{eqnarray*}
\end{proof}

  \begin{proof}[Proof of Lemma \ref{lem:templem}]

  We consider two cases :
  (i) $4\lambda_n \sum_{j \in \mathcal{M}_*}  d^{1/\gamma_j^*} \rho(\| {\boldsymbol\Delta}_j\|_{\gamma_j}) \geq  \| f_{ \mathbf b} - f_{\boldsymbol \beta^*} \|_{n,\mathbf b^*}^2$,
 and (ii) $4\lambda_n \sum_{j \in \mathcal{M}_*}   d^{1/\gamma_j^*} $ $ \rho(\| {\boldsymbol\Delta}_j\|_{\gamma_j}) \leq \| f_{ \mathbf b} - f_{\boldsymbol \beta^*} \|_{n,\mathbf b^*}^2$.

 {\bf Case (i)} From \eqref{eq:temp19}, we have
$$
 \| f_{\widehat{{\boldsymbol \beta}}} - f_{{\boldsymbol \beta}^*}\|_{n,\mathbf b_{\widehat { \boldsymbol\beta}} }^2  + \lambda_n \sum_{j=1}^p   d^{1/\gamma_j^*}\rho(\| {\boldsymbol\Delta}_j \|_{\gamma_j}) \leq  8 \lambda_n \sum_{j \in \mathcal{M}_*}  d^{1/\gamma_j^*} \rho(\| {\boldsymbol\Delta}_j\|_{\gamma_j}).
$$
  This implies that $ \sum_{j \in \mathcal{M}_*^c}  d^{1/\gamma_j^*} \rho(\| {\boldsymbol\Delta}_j\|_{\gamma_j}) < 7  \sum_{j \in \mathcal{M}_*}  d^{1/\gamma_j^*} \rho(\| {\boldsymbol\Delta}_j\|_{\gamma_j}) $ or that ${\boldsymbol\Delta} =\hat \bbeta - \mathbf b \in \mathbb{C}_{7,\rho}$ as defined in RE condition. For such ${\boldsymbol\Delta}$,    from the RE condition in \eqref{eq:RE} we have with
  $$\bar d=\sum_{j \in \mathcal{M}_*} d^{2/\gamma_j^*}, $$
    \begin{eqnarray}\nonumber
 {\| f_{\hat {\boldsymbol\beta} } - f_{ {\boldsymbol\beta}^*}\|_{n,\mathbf b_{\widehat { \boldsymbol\beta}} }^2 \leq
 8 \lambda_n \sqrt{\bar d} \sqrt{ \sum_{j \in \mathcal{M}_*} \rho^2(\| {\boldsymbol\Delta}_j\|_{\gamma_j}) }
  \leq  8 \lambda_n \frac{\sqrt{\bar d} }{\zeta} \sqrt{\bDelta^T \bigtriangledown^2 \mathcal{L}_n(\bbeta^*) \bDelta}}
  \end{eqnarray}
  The left hand side can be further bounded  using   Proposition 2 and triangle inequality with
     \begin{eqnarray}\nonumber
        8 \lambda_n    \frac{\sqrt{\bar d} }{\zeta} \left(  {  \exp\{  a_{\hat{\boldsymbol\beta} - \boldsymbol\beta^*}\} \| f_{{\boldsymbol\beta}^*} - f_{\hat {\boldsymbol\beta} }\|_{n,\mathbf b_{\widehat { \boldsymbol\beta}} }  +   \exp\{ a_{\mathbf b -\boldsymbol\beta^*}  \}  \| f_{\mathbf b} - f_{  {\boldsymbol\beta}^* }\|_{n, \mathbf b^*}   }  \right)
         \end{eqnarray}
Furthermore, with  the simple inequality  $a b \leq b^2/2 + a^2/2$, we can further upper bound the left hand side with
   \begin{eqnarray}\nonumber
      \leq
    16 \lambda_n^2   \frac{ {\bar d} }{\zeta^2}  \exp\{ 2a_{\mathbf b -\boldsymbol\beta^*}\} + \| f_{\mathbf b} - f_{  {\boldsymbol\beta}^* }\|_{n, \mathbf b^*}^2   +  32\lambda_n^2   \frac{ {\bar d} }{\zeta^2}  \exp\{ 2 a_{\hat{\boldsymbol\beta} - \boldsymbol\beta^*}\}  + \frac{1}{2}\| f_{{\boldsymbol\beta}^*} - f_{\hat {\boldsymbol\beta} }\|_{n,\mathbf b_{\widehat { \boldsymbol\beta}} }^2  .
     \end{eqnarray}
 Combining all of the above we obtain
        \begin{eqnarray*}
  {\| f_{\hat {\boldsymbol\beta} } - f_{ {\boldsymbol\beta}^*}\|_{n,\mathbf b_{\widehat { \boldsymbol\beta}} }^2 \leq    2\| f_{\mathbf b} - f_{  {\boldsymbol\beta}^* }\|_{n, \mathbf b^*}^2   +  64\lambda_n^2   \frac{ {\bar d} }{\zeta^2}  \exp\{ 2 a_{\hat{\boldsymbol\beta} - \boldsymbol\beta^*}\}    +32 \lambda_n^2   \frac{ {\bar d} }{\zeta^2}  \exp\{ 2a_{\mathbf b -\boldsymbol\beta^*}\}.}
  \end{eqnarray*}

 To upper bound the LHS of the previous inequality we  bound the two exponential terms independently. .

 First: let $b= \bbeta^*$ in \eqref{eq:temp19}. Then, by using the RE condition  and all equations above we obtain that $y_1=\sum_{j \in \mathcal{M}^*} d^{1/\gamma_j^*}\rho(\|\hat \bbeta_j - \bbeta^*_j \|_{\gamma_j})  \geq 0$ and $\upsilon_1=\sum_{j \in \mathcal{M}^*} \|\hat \bbeta_j - \bbeta^*_j \|_{\gamma_j}  \geq 0 $
are such that $  a_{\hat{\boldsymbol\beta} - \boldsymbol\beta^*} \leq C \upsilon_1$  and  \[
y_1\exp\{-2 C \upsilon_1\}  \leq 16 \lambda_n^2   \frac{ {\bar d} }{\zeta^2}.
 \]
 From convexity of $\rho$ we know that $\rho(\|\hat \bbeta_j - \bbeta^*_j \|_{\gamma_j}) \geq  \rho'(0+) \|\hat \bbeta_j - \bbeta^*_j \|_{\gamma_j} $, hence $\upsilon_1,\upsilon_2$
  satisfy  $\upsilon_1 \geq  \rho'(0+)    \upsilon_2   $. Combining all  the above,  $\upsilon_1$ solves
 \begin{equation}\label{eq:28}
\upsilon_1\exp\{-2 C \upsilon_1\}  \leq 16 \lambda_n^2   \rho'(0+)\frac{ {\bar d} }{\zeta^2}.
 \end{equation}

 Second, we consider the case of general $\mathbf b $ possibly different from $\bbeta^*$.
 In such cases,
  $$ \| f_{ \mathbf b} - f_{\boldsymbol\beta^*} \|_{n,\mathbf b^*}^2 \leq 4\lambda_n \sum_{j \in \mathcal{M}_*}  d^{1/\gamma_j^*} \rho(\| \hat{\boldsymbol\beta}_j - \bbeta^*_j\|_{\gamma_j}) +4\lambda_n \sum_{j \in \mathcal{M}_*}  d^{1/\gamma_j^*} \rho(\| \mathbf b_j - \bbeta^*_j\|_{\gamma_j})$$
  Then, by utilizing Proposition 2 on the left and Caushy-Shwarz inequality to the right, we notice that
  $$  \exp\{ -2a_{\mathbf b -\boldsymbol\beta^*}\} \zeta^2  \sum_{j \in \mathcal{M}_*}   \rho^2(\| \mathbf b_j - \bbeta^*_j\|_{\gamma_j}) \leq \| f_{ \mathbf b} - f_{\boldsymbol\beta^*} \|_{n,\mathbf b^*}^2 \leq 4 \lambda_n  \bar d  \sqrt{ \upsilon_2 } +4\lambda_n  \bar d \sqrt{\sum_{j \in \mathcal{M}_*}    \rho^2(\| \mathbf b_j - \bbeta^*_j\|_{\gamma_j}) }.$$
To that end, let us denote with $y_2= \sum_{j \in \mathcal{M}_*}   \rho^2(\| \mathbf b_j - \bbeta^*_j\|_{\gamma_j}) \geq 0$ and  $\upsilon_2= \sum_{j \in \mathcal{M}_*}    \| \mathbf b_j - \bbeta^*_j\|_{\gamma_j}^2 \geq 0$ and observe that  $a_{\mathbf b -\boldsymbol\beta^*} \leq C  \upsilon_2$
\[
\zeta^2   y_2\exp\{ -2 C \upsilon_2\}  - 4\lambda_n  \bar d  \sqrt{y_2} \leq 4 \lambda_n  \bar d  \sqrt{ \upsilon_1 } .
\]
Utilizing the equation $\upsilon_1$ satisfies and the convexity of $\rho$ we have
\begin{equation}\label{eq:29}
 \upsilon_2  \exp\{ -2 C \upsilon_2\}  - 4\lambda_n  \frac{\bar d }{ \zeta^{2} {\rho'}^{2}(0+)} \sqrt{\upsilon_2} \leq 16 \lambda_n^2  \frac{{\bar d}^{3/2}}{    {\rho'}^{3/2}(0+) \zeta^{3}}.
\end{equation}
Although $\upsilon_2$ depends on $\bb$, we observe that the previous inequality holds uniformly over $\bb$ hence we have suppressed the dependence on $\bb$ in the notation of $\upsilon_2$.

  {\bf Case (ii)} From \eqref{eq:temp19}, we have
$$
 \| f_{\widehat{{\boldsymbol \beta}}} - f_{{\boldsymbol \beta}^*}\|_{n,\mathbf b_{\widehat { \boldsymbol\beta}} }^2    \leq  2 \| f_{ \mathbf b} - f_{\boldsymbol \beta^*} \|_{n,\mathbf b^*}^2
 \leq
  64\lambda_n^2   \frac{ {\bar d} }{\zeta^2}  \exp\{ 2 C \upsilon_1\}    +32 \lambda_n^2   \frac{ {\bar d} }{\zeta^2}  \exp\{ 2 C \upsilon_2 \}  + 2  \min_{\mathbf b \in \mathcal{B}}  \| g- f_{\mathbf b}\| ^2  .
$$

\end{proof}

\begin{proof}[Proof of Lemma \ref{lem:globalSOIb}]

     Following the same steps as in the proof of Lemma \ref{lem:templem},   we obtain easily that
     ${\boldsymbol\Delta} \in \mathbb{C}_3$ for ${\boldsymbol\Delta}=\widehat{\bbeta}-\bbeta^*$ (exact steps are omitted).
Combined with assumption RE($7,s,\bgamma$)  \eqref{eq:RE}, it leads to
   \begin{equation}\label{eq:SOI2a}
  \| f_{\widehat{{\boldsymbol \beta}}} - f_{{\boldsymbol \beta^*}}\|_{n,\widehat{\boldsymbol \beta }^*}^2 \leq 32 \frac{  \lambda_n ^2}{ \zeta^2} e^{2 \upsilon_1}\sum_{j\in \mathcal{M}_*} d^{2/\gamma_j^*} ,
  \end{equation}
  for $0 \leq \upsilon_1 \leq 1$  satisfying \eqref{eq:upsilons}.
  This result gives a preliminary step towards the final statement. The right- hand side is a complicated random norm  (introduced in (11)). The rest of the proof establishes tight non-trivial lower bounds on its size.
Together with Proposition \ref{lem:approx}, we have
  \begin{eqnarray} \nonumber
32\frac{ \lambda_n^2 }{\zeta^2} e^{2 \upsilon_1}\sum_{j\in \mathcal{M}_*} d^{2/\gamma_j^*}   &\geq & \| f_{\widehat{{\boldsymbol\beta}}} - f_{{\boldsymbol \beta}^*}\|_{n,\widehat{\boldsymbol \beta }^*}^2 = \frac{ \int_0^\tau (\widehat {\boldsymbol \beta}- \boldsymbol{\beta}^*)^T \bV_n(\bb_{\widehat{\boldsymbol\beta}},t)  (\widehat {\boldsymbol \beta} - \boldsymbol{\beta}^*)d\bar{N}(t) }{\|\widehat {\boldsymbol \beta}_{\mathcal{M}_*} - \boldsymbol{\beta}^*_{\mathcal{M}_*} \|_{1,\boldsymbol\gamma}^2} \|\widehat {\boldsymbol \beta}_{\mathcal{M}_*} - \boldsymbol{\beta}^*_{\mathcal{M}_*} \|_{1,\boldsymbol\gamma}^2 \\
 &\geq& e^{-2a_{\hat {\boldsymbol\beta} - \boldsymbol\beta^*}} \zeta^2 \|\widehat {\boldsymbol \beta}_{\mathcal{M}_*} - \boldsymbol{\beta}^*_{\mathcal{M}_*} \|_{1,\boldsymbol\gamma}^2, \label{eq:34}
 \end{eqnarray}
   where we used the notation $  \|\widehat {\boldsymbol \beta} _{\mathcal{M}_*} - \boldsymbol{\beta}^*_{\mathcal{M}_*}  \|_{1,\boldsymbol \gamma}^2 =\sum_{j \in{\mathcal{M}_*} } \| \widehat {\boldsymbol \beta}_j - \boldsymbol{\beta}^*_j \|_{\gamma_j}^2$, and $a_{\hat {\boldsymbol\beta} - \boldsymbol\beta^*}=\max_{1 \leq q, i \leq n} | (\hat {\boldsymbol\beta} - \boldsymbol\beta^*)^T \Psi(\bX_i) - \Psi(\bX_q)|\leq 2 C \| \hat{\boldsymbol \beta}  - \boldsymbol{\beta}^* \|_1 $. The rest of the proof is based on the analysis of the upper bound for the norm $\| \hat{\boldsymbol \beta}  - \boldsymbol{\beta}^* \|_1 $. The goal is to first  find the worst case upper bound that satisfies \eqref{eq:34}. Therefore,  the  desired upper bound is  the optimal solution of  the following optimization problem
      \[
      \begin{array}{cc}
      \max & \| \bx\|_1
      \ \ \ \\
      & \\
      \mbox{s.t.} & \ \ \  e^{- \| \bx\|_1}\| \bx\|_{1,\gamma} \leq z,
      \end{array}
      \]
      for $z= 16\frac{ \lambda_n^2 }{\zeta^4} e^{2 \upsilon_1} \sum_{j\in \mathcal{M}_*} d^{2/\gamma_j^*}  $. Because
   $
   \| \bx\|_{1,\gamma}  \geq  d^{- 1} \|\bx \|_{1},
   $
   the optimal value of the previous problem is  upper bounded by the optimal value of the following problem
   \[
    \begin{array}{cc}
      \max & u \  \ \ \ \\
      & \\
      \mbox{s.t.} &\ \ \  e^{- u}u \leq  zd , \\
      & \ \ \ u \geq 0.
      \end{array}
      \]
Function $e^{- u}u$  is neither convex or concave.   It is concave up to $u=2$ and then convex with exponentially rate of convergence towards zero. When $zd>1/e$, the optima is reached at $u=1$.  When $zd <2 e^{-2}$, $u \to \infty$ exponentially fast.
Thus, for $\lambda_n$ satisfying
\[
 e^{-1} \zeta^4 \leq 32 { \lambda_n^2 }{ e^{2 \upsilon_1}}d \sum_{j\in \mathcal{M}_*} d^{2/\gamma_j^*},
\]
we have $\| \hat{\boldsymbol \beta}  - \boldsymbol{\beta}^* \|_1  \leq 1$.  Under such conditions for some constant $c_0 > 1$,
$
a_{\hat {\boldsymbol\beta} - \boldsymbol\beta^*} \leq 2 C,
$ and we utilize \eqref{eq:34} to conclude
   \[
    \|\widehat {\boldsymbol \beta}_{\mathcal{M}_*} - \boldsymbol{\beta}^*_{\mathcal{M}_*} \|_{1,\boldsymbol\gamma}^2 \leq 32e^{2 C + 2\upsilon_1} \frac{ \lambda_n^2 }{\zeta^4}\sum_{j\in \mathcal{M}_*} d^{2/\gamma_j^*}.
   \]
   From Cauchy Schwartz inequality, we have that
$
     \sum_{j\in \mathcal{M}_*}  d^{1/\gamma_j^*}\| \widehat {\boldsymbol \beta}_j - \boldsymbol{\beta}^*_j \|_{\gamma_j} $ is less than or equal to $\sqrt{\sum_{j \in \mathcal{M}_*} d^{2/\gamma_j^*}}
     \sqrt{ \sum_{j\in \mathcal{M}_*} \| \widehat {\boldsymbol \beta}_j - \boldsymbol{\beta}^*_j \|_{\gamma_j}^2 }$ which is according to inequality above   upper bounded with
     $ 4e^C \frac{ \lambda_n }{\zeta^2} \sum_{j\in \mathcal{M}_*} d^{2/\gamma_j^*}
$.
Knowing that $\widehat \bbeta - \bbeta^* \in \mathbb C_3$ and using the convexity of $\rho$, we have $\| \rho(\widehat {\boldsymbol \beta} _{\mathcal{M}_*^c} - \boldsymbol{\beta}^*_{\mathcal{M}_*^c})\|_{1} \leq 3 \|\rho(\widehat {\boldsymbol \beta} _{\mathcal{M}_*} - \boldsymbol{\beta}^*_{\mathcal{M}_*})\|_{1}$ and  thus
   \begin{eqnarray*}
      \sum_{j=1}^p d^{1/\gamma_j^*} \| \widehat {\boldsymbol \beta}_j - \boldsymbol{\beta}^*_j \|_{\gamma_j}
      &\leq& 16 \sqrt{2} e^{C+\upsilon_1} \frac{ \lambda_n }{  \zeta^2} \sum_{j\in \mathcal{M}_*} d^{2/\gamma_j^*} .
      \end{eqnarray*}
        \end{proof}

 \begin{proof}[Proof of Lemma \ref{lemma:smooth}]
 Let

 $$\mathcal{T}_n = \{\| \tilde{\bh}_{n,j}(\bbeta^*)\|_{\gamma_j^*} \leq  2 \lambda_n  \max\{d^{1/\gamma_j^*} \sqrt{d}\} \min_{1 \leq j \leq p}  \lambda_{\min}^{}(\bR_j)  \rho'(0+), \forall j \in \{1,\cdots, p\} \},$$
 with
$\tilde{\bh}_{n,j}(\bbeta^*)= - \frac{1}{n} \sum_{i=1}^n \int_0^\tau ({\widetilde \bE}_{n,j}(\bbeta^*,t) -\bR_j^{-1} \bPsi(X_{ij}))dM_i(t)$,
 \[
 {\widetilde \bE}_{n,j}(\bbeta^*,t) =\frac{1}{n} \sum_{i=1}^n  \frac{Y_i(t) \bR_j^{-1} \bPsi(X_{ij}) }{\frac{1}{n} \sum_{l=1}^n Y_l(t) \exp\{ \sum_{j=1}^p  {\bbeta^*_j}^T \bR_j^{-1} \bPsi(X_{lj})\}} \exp\{ \sum_{j=1}^p  {\bbeta^*_j}^T \bR_j^{-1} \bPsi(X_{ij})\}
 \]

 We first adapt the results of  Lemma \ref{lemma:min1} with the following few steps
 \begin{eqnarray*}
 - {\bb}_j^T \bh_{n,j}(\bbeta^*)+\lambda_n  \sqrt{d}  \rho\left(\|\bb_j \|_{\gamma_j} +\|\bb_j \|_{2} \right) \\
 \geq \|\bb_j\| _{\gamma_j} \left( -\|\bh_{n,j}(\bbeta^*)\|_{\gamma_j^*}+\lambda_n  \sqrt{d}  \rho'(0+) \left( 1 + \frac{\|\bb_j  \|_2}{\|\bb_j\| _{\gamma_j}} \right) \right).
 \end{eqnarray*}
 For $\gamma_j \geq 2 $, we know that $\|\bb_j\| _{\gamma_j} \leq \|\bb_j  \|_2$. This relation leads to the conclusion that previous quantity is lower bounded with
 \[
 \geq \|\bb_j\| _{\gamma_j} \left( -\|\bh_{n,j}(\bbeta^*)\|_{\gamma_j^*}+ 2 \lambda_n  \sqrt{d}  \rho'(0+)  \right),
 \]
 which leads us to conclude that the results of Lemma \ref{lemma:min1} hold for this particular penalty.  Size of the set $\mathcal{T}_n $   is  easily deducible by adapting the very last proof of Theorem \ref{cor:localSOI} (exact details are omitted).

To prove equivalent results to those of Section \ref{sec:lan}, we need to define new constants corresponding to  $a_{\mathbf{v}}$ and
 $\underline{\underline{\omega}}$.
First, the equivalent of $\bV_n(\bb)$ has  extra  $\bR_j^{-1}$ terms, which will factor into f  $a_i$ terms  (of Proposition \ref{lem:approx}) as $(\bb-\bbeta^*)( \bR^{-1}\bPsi(\bX_i)- \bE_n(\bbeta^*,t))$. $\bR$ is a diagonal block matrix
 $$\bR= \left( \begin{array}{cccc} \bR_1& \mathbf 0& \cdots& \mathbf 0 \\ \mathbf 0 & \bR_2 & \cdots & \mathbf 0 \\ \vdots \\ \mathbf 0&  \mathbf 0 & \cdots& \bR_p \end{array}\right).$$
Second, as $\bar{\mathbf v}_j = \mathbf v_j \bR_j^{-1}$
 \[
 \bar a_{\mathbf{v}} =\max_{1 \leq i,q \leq n} \biggl | {\bar{\mathbf v}}^T (\bPsi(\bX_i) - \bPsi(\bX_q))  \biggl| \leq \max_{1 \leq j \leq p} r(\bR_j^{-1}) a_{\mathbf v}
 \]
 with spectral radius $r(\bR_j^{-1}) = \max_{k=1,\cdots,d}|\lambda_k(\bR_j^{-1})| =  \max_{k=1,\cdots,d}|\lambda_k(\bR_j)|^{-1}= \lambda_{\min}^{-1}(\bR_j)$. Then, $ \bar a_{\mathbf{v}}  \leq \max_{1 \leq j \leq p}  \lambda_{\min}^{-1}(\bR_j) a_{\mathbf v}$. Thus, the result of Proposition  \ref{lem:approx} follows with $\eta$ equal to $\max_{1 \leq j \leq p}  \lambda_{\min}^{-1}(\bR_j) a_{\mathbf v}$.

 The definition of the weights, $\omega_i(\bb)$, in the proof of Proposition \ref{prop:approx} will be changed to address the new weighting matrix, $\bR_j$,. Once they are redefined
 with

 $$\underline{\underline \omega}_{\mathbb{S}}:=\min_{i \in \{1,\cdots,n\}, i \in \cup_{q=1}^n \mathbb{R}_q} \left\{ \frac{\sum_{q=1}^N \exp\{ \sum_{j=1}^p {\bbeta^*_j}^T \bR_j^{-1} \bPsi(X_{ij})\} 1\{i \in \mathbb{R}_q\}\}}{\sum_{l \in \mathbb{R}_q} \exp\{ \sum_{j=1}^p {\bbeta^*_j}^T \bR_j^{-1} \bPsi(X_{lj})\}} \right\},$$
 the exact steps of the proof   of Proposition \ref{prop:approx}   will follow easily, and thus we omit the details here.

   \end{proof}

\end{document}